\documentclass[11pt, a4paper, reqno]{article}

\RequirePackage{amsthm,amsmath,amsfonts,amssymb}
\RequirePackage[numbers]{natbib}
\RequirePackage[colorlinks,citecolor=blue,urlcolor=blue]{hyperref}
\RequirePackage{graphicx}

\theoremstyle{plain}
\newtheorem{theo}{Theorem}[section]
\newtheorem{prop}[theo]{Proposition}
\newtheorem{coro}[theo]{Corollary}
\newtheorem{lemme}[theo]{Lemma}

\theoremstyle{remark}
\newtheorem{definition}[theo]{Definition}

\newtheorem{rem}[theo]{Remark}

\newtheorem{ass}[theo]{Assumption}

\usepackage{comment}
\usepackage{dsfont}
\usepackage{enumitem}

\newcommand{\w}{\widehat}
\newcommand{\wt}{\widetilde}
\newcommand{\one}{\mathds{1}}

\newcommand{\E}{\mathbb{E}}

\newcommand{\R}{\mathbb{R}}

\newcommand{\Var}{\text{Var}}
\newcommand{\Cov}{\text{Cov}}

\newcommand{\dd}{\mathrm{d}}
\newcommand{\Tm}{T_{\min}}
\newcommand{\norm}[1]{\left\Vert #1\right\Vert}

\title{Nonparametric estimation of the stationary density for Hawkes-diffusion systems with known and unknown intensity}

\author{Chiara Amorino \thanks{Universitat Pompeu Fabra and Barcelona School of Economics, Department of Economics and Business, Ram\'on Trias Fargas 25-27, 08005, Barcelona, Spain}, Charlotte Dion-Blanc \thanks{LPSM, Sorbonne Université, Université Paris Cité, CNRS, F-75005 Paris, France.}, Arnaud Gloter \thanks{LaMME, CNRS, Univ Evry, Universit\'e Paris-Saclay, 91037, Evry, France}, 
Sarah Lemler \thanks{MICS, \'Ecole CentraleSup\'elec, Universit\'e Paris-Saclay, France }}

\begin{document}

\maketitle

\begin{abstract}
We investigate the nonparametric estimation problem of the density $\pi$, representing the stationary distribution of a two-dimensional system $\left(Z_t\right)_{t \in[0, T]}=\left(X_t, \lambda_t\right)_{t \in[0, T]}$. In this system, $X$ is a Hawkes-diffusion process, and $\lambda$ denotes the stochastic intensity of the Hawkes process driving the jumps of $X$. Based on the continuous observation of a path of $(X_t)$ over $[0, T]$, and initially assuming that $\lambda$ is known, we establish the convergence rate of a kernel estimator $\widehat\pi\left(x^*, y^*\right)$ of $\pi\left(x^*,y^*\right)$ as $T \rightarrow \infty$. Interestingly, this rate depends on the value of $y^*$ influenced by the baseline parameter of the Hawkes intensity process. From the rate of convergence of $\widehat\pi\left(x^*,y^*\right)$, we derive the rate of convergence for an estimator of the invariant density $\lambda$. Subsequently, we extend the study to the case where $\lambda$ is unknown, plugging an estimator of $\lambda$ in the kernel estimator and deducing new rates of convergence for the obtained estimator. 
The proofs establishing these convergence rates rely on probabilistic results that may hold independent interest.
We introduce a Girsanov change of measure to transform the Hawkes process with intensity $\lambda$ into a Poisson process with constant intensity. To achieve this, we extend a bound for the exponential moments for the Hawkes process, originally established in the stationary case, to the non-stationary case. 
Lastly, we conduct a numerical study to illustrate the obtained rates of convergence of our estimators.
\end{abstract}

\begin{abstract}[language=french]
 Nous nous intéressons au problème d'estimation non paramétrique de la densité $\pi$, qui représente la distribution stationnaire d'un système bidimensionnel $\left(Z_t\right)_{t \in[0, T]}=\left(X_t, \lambda_t\right)_{t \in[0, T]}$. Dans ce système, $X$ est un processus de diffusion à sauts, et $\lambda$ désigne l'intensité stochastique du processus de Hawkes qui modélise les sauts de $X$. En nous basant sur l'observation continue d'une trajectoire de $(X_t)$ sur $[0, T]$, et en supposant initialement que $\lambda$ est connue, nous étudions la vitesse de convergence d'un estimateur à noyau $\widehat\pi\left(x^*, y^*\right)$ de $\pi\left(x^*,y^*\right)$ lorsque $T \rightarrow \infty$. Nous avons mis en évidence que ce taux dépend de la valeur prise par $y^*$ en fonction du paramètre de base de l'intensité du processus de Hawkes. À partir de la vitesse de convergence de $\widehat\pi\left(x^*,y^*\right)$, nous déduisons la vitesse de convergence d'un estimateur de la densité invariante $\lambda$. Par la suite, nous étendons l'étude au cas où $\lambda$ est inconnue, en introduisant un estimateur de $\lambda$ dans l'estimateur à noyau et en déduisant de nouvelles vitesses de convergence pour l'estimateur obtenu.
Les démonstrations établissant ces vitesses de convergence reposent sur des résultats probabilistes qui présentent un intérêt en tant que tels. Nous introduisons un changement de mesure de Girsanov pour transformer le processus de Hawkes d'intensité $\lambda$ en un processus de Poisson d'intensité constante. Pour ce faire, nous étendons une borne sur les moments exponentiels pour le processus de Hawkes, initialement établie dans le cas stationnaire, au cas non stationnaire.
Enfin, nous réalisons une étude numérique pour illustrer les vitesses de convergence obtenues de nos estimateurs.
\end{abstract}



\section{Introduction}

In this work, we study a {time-homogeneous} Hawkes-diffusion process, where the dynamics of the process \((X_t)\) are governed by the stochastic differential equation:  
\begin{equation}\label{Eq : model diff X intro}
\dd X_t = b(X_t)\, dt + \sigma(X_t)\, \dd W_t + a(X_{t-})\, \dd N_t,
\end{equation}
with initial condition \(X_0\). Here, \(W\) represents a standard Brownian motion independent of the jump process \(N\), which is a one-dimensional linear exponential Hawkes process with intensity \((\lambda_t)\), defined later. For a comprehensive review of this process, we refer to \cite{jumpdiffreview}.  

The Hawkes process, introduced by Hawkes in \cite{hawkes1971spectra}, has been extensively studied in various contexts, with applications spanning fields such as seismology \cite{ogata1988statistical}, sociology \cite{crane2008robust}, and finance \cite{bacry2015hawkes}. Notably, it appears in financial modeling in \cite{ait2015modeling} and has been recently explored in \cite{ait2023saddlepoint} for risk management. In \cite{mancini}, the dynamics of log-returns on the S\&P500 index were modeled within this framework, while a stochastic control perspective was examined in \cite{bensoussan2024stochastic}.  

This process is particularly relevant for modeling continuous phenomena influenced by external events with an auto-excitation structure. For instance, it is applicable in interest rate modeling in insurance, as shown in \cite{gomez2016estimation}. However, its most prominent applications lie in neuroscience, where it is used to model the evolution of the membrane potential influenced by signals from surrounding neurons (see \cite{DL} and \cite{tchouanti2024separation} for details). Specifically, it is common to describe a neuron's spike train using a Hawkes process to capture the auto-excitation: for certain neuron types, one spike increases the probability of subsequent spikes.  

Moreover, this process generalizes Poisson jumps (or L\'evy jumps, which have independent increments) by incorporating auto-exciting jumps. Consequently, many practical applications of L\'evy processes can also serve as motivations for this study. Additionally, the Hawkes-diffusion process offers a more tractable alternative compared to models driven by L\'evy processes, further broadening its applicability.

In this paper, we focus on the pointwise estimation of the invariant density $\pi(x^*,y^*)$ of the couple $(X,\lambda)$ using kernel density estimators.
Kernel density estimation is a widely used estimation method for stochastic processes. For instance, it has been applied to estimate the marginal density of continuous-time processes \cite{banon1978nonparametric,bosq1998minimax} within jump-diffusion frameworks \cite{amorino2021invariant, amorino2022optimal,mancini2011threshold}, and in models driven by fractional Brownian motion \cite{amorino2024fast, bertin2020adaptive}. A rich body of literature explores the convergence rates of kernel estimators for invariant densities of various diffusion processes, motivating our approach to estimating the invariant density $\pi$ for the pair \((X, \lambda)\).

\paragraph{Related works.}
The literature on statistical inference for Hawkes diffusion processes is extensive. Hawkes processes are defined by their intensity function, which describes the infinitesimal probability of a new event occurring, conditional on the past. Consequently, estimating these processes primarily involves estimating this function. A variety of techniques exists, depending on the model’s flexibility and the available data.  
In the parametric setting, the standard approach relies on maximum likelihood estimation (see \cite{bonnet2023inference, clinet2017statistical, Ogata}). Other methods, such as the EM algorithm \cite{veen2008estimation} or least squares \cite{bacry2015hawkes, DDBLS}, are also commonly employed. On the other hand, when the intensity is defined in a non-parametric way, common approaches include Bayesian inference \cite{deutsch2022estimating, sulem2024bayesian}, non-parametric EM algorithms \cite{lewis2011nonparametric}, and RKHS techniques \cite{yang2017online}. Furthermore, parameter estimation for marked Hawkes processes has been explored in works such as \cite{bonnet2024testing, clinet2022quasi, goda2023sparse}. We direct interested readers to these papers and the references therein for further details.  

Nonparametric estimation of the coefficients for a Hawkes-diffusion process  can be found in \cite{ADBGL} and \cite{DL}. These works proposed nonparametric estimators for the coefficients \(b\), \(\sigma\), and \(a\), with rigorous error bounds.  

Despite the extensive body of work cited above, we are unaware of any studies addressing the estimation of the invariant density of the pair \((X, \lambda)\) in our setting. This paper aims to fill this gap. Nevertheless, invariant density estimation is a statistical problem that has been extensively studied in other frameworks.

In the one-dimensional setting, it is well established in the field of statistics for stochastic processes that the estimation of the invariant density can be achieved under standard nonparametric assumptions with a parametric rate (see Chapter 4.2 in \cite{kutoyants2013statistical} for further details, as well as \cite{comte2002adaptive, comte2005super}, while \cite{schmisser2013nonparametric} studies the estimation of derivatives of the invariant density).  

For multidimensional diffusion processes, a key reference is the work \cite{dalalyan2007asymptotic}, which demonstrated that the convergence rate for estimating the invariant density belonging to an isotropic H\"older class is \(T^{-\frac{2 \beta}{2 \beta + d - 2}}\), where $\beta$ is the regularity. These results have been extended to the anisotropic setting in \cite{strauch2018adaptive}, with the rate proven to be optimal in a minimax sense in \cite{amorino2024minimax}. Similar developments for stochastic differential equations driven by other types of noise include \cite{amorino2021invariant, amorino2022optimal} for SDEs driven by Lévy processes, and \cite{amorino2024fast, bertin2020adaptive} for SDEs driven by fractional Brownian motion.  

The paper \cite{delattre2022rate} investigates the estimation of the invariant density for a two-dimensional hypoelliptic stochastic system. Notably, this work shares several similarities with ours, as discussed further below Corollary \ref{coro:rates}. This resemblance arises because our system also exhibits hypoelliptic features: the stochastic differential equation governing the intensity of the Hawkes process is degenerate, lacking a Brownian component.  

Lastly, we highlight the work \cite{krell2021nonparametric}, which examines the nonparametric estimation of jump rates for a specific class of one-dimensional piecewise deterministic Markov processes (PDMPs) and constructs an estimator for the stationary density. As Section \ref{sec:estLambda} of this paper focuses on estimating the invariant density of the Hawkes process intensity \(\lambda\), a natural comparison between the results in \cite{krell2021nonparametric} and our findings for our particular PDMP process \(\lambda\) arises. This comparison is elaborated upon in the discussion following Corollary \ref{c: rate lambda}.


\paragraph{Our contribution.}
We investigate two observation schemes. In the first, we assume that the intensity process \((\lambda_t)\) is known. Utilizing continuous observation of \((X_t, N_t)_{t \in [0, T]}\), we propose a kernel estimator and analyze its convergence rate as the time horizon \(T\) tends to infinity. Since \(N\) is the jump process of \(X\), the events of \(N\) are entirely determined by the continuous observation of \(X\).

We demonstrate that the convergence rates of the kernel estimator depend on the value of \(y\) relative to the baseline parameter \(\xi\) of the Hawkes process. Specifically, as is standard in the literature, we assume that the invariant density belongs to an anisotropic H\"older class with regularity \(\beta = (\beta_1, \beta_2)\) (see, for example, \cite{amorino2021invariant, hopfner2002non, juditsky2002nonparametric, lepski1997optimal, tsybakov1997nonparametric}). After selecting the rate-optimal bandwidths \(h_1\) and \(h_2\), we prove that, for any arbitrarily small \(\varepsilon > 0\) for  \(y^* \in[\xi,+\infty)\), we have the following:

\[
\E\left[(\w{\pi}_{h_1, h_2}(x^*,y^*)-{\pi}(x^*,y^*))^2\right] \lesssim 
\begin{cases}
\left(\frac{1}{T}\right)^{\frac{2 \beta_2}{2 \beta_2 + 1}} \log(T), &\quad y^* = \xi, \\
\left(\frac{1}{T}\right)^{\frac{4 \beta_2}{4 \beta_2 + 1} - \varepsilon}, &\quad y^* > \xi.
\end{cases}
\]

Interestingly, the estimation procedure performs better for \(y^* > \xi\). The intuition behind this result is that the deterioration in the estimation procedure for \(y^* = \xi\) occurs because the process rarely visits the baseline (see the discussion following Remark \ref{rk: idea proof girsanov} for more details).

The convergence rates established above are both surprising and fast. We believe that the bounds we have proven are sharp, even though we have not yet been able to rigorously address the optimality in a minimax sense by establishing a corresponding lower bound. A detailed explanation of the challenges associated with deriving a lower bound is provided in Section \ref{s: comments}.

Below Corollary \ref{coro:rates}, we provide a thorough comparison of our results with the convergence rates for invariant density estimation in related literature. This comparison highlights how our results effectively separate the contributions of the jump intensity \(\lambda\) from those of the diffusion \(X\). It underscores the fact that the rates we derive capture an intrinsic characteristic of the SDE driven by the Hawkes process.
 
Additionally, we propose an estimator for the invariant density of the intensity \((\lambda_t)\) using a kernel density estimator based solely on the observation of the process \(\lambda\). We derive its mean squared error convergence rate. As before, the convergence rate depends on whether \(y\) exceeds the baseline \(\xi\) or not. After selecting the rate-optimal bandwidth \(h_2\), the rate is given by: 
\[
\E[(\hat{\pi}_{h_2}(y^*) - \pi(y))^2] \lesssim 
\begin{cases}
 \left(\frac{1}{T}\right)^{\frac{2 \beta_2}{2 \beta_2 + 1}} \log(T), & \text{if } y^* = \xi, \\
\frac{1}{T}(\log T + T^{\varepsilon}), & \text{if } y^* > \xi.
\end{cases}
\]
The first noteworthy observation is that, as before, the estimation procedure performs better when \(y^* > \xi\). In particular, in this case, we recover, up to logarithmic terms, a superoptimal convergence rate. Another important remark concerns the comparison with \cite{krell2021nonparametric}, where the authors, studying the estimation of the invariant density of one-dimensional PDMPs, achieved a classical nonparametric convergence rate. Interestingly, this rate aligns with the one we obtain for \(y^* = \xi\). However, in our case, this rate corresponds only to the worst-case scenario.

In the second observation scheme, we consider a more general setting where the intensity process \((\lambda_t)\) is unknown and must be estimated. Based on the observed jump times of \(N\) (and consequently \(X\)), we estimate the parameters governing the intensity process \((\lambda_t)\), ultimately deriving an estimator for \(\lambda\) at any given time. In Theorem \ref{th: bound lambda hat}, we establish that this estimator achieves a parametric rate, a result we believe may hold independent interest for future studies.

Using this estimator, we construct a plug-in kernel estimator for \(\pi(x^*, y^*)\) and determine its convergence rates. Unlike the previous scenario, these rates are unaffected by whether or not \(y^* = \xi\). However, as one may expect, the rates are slower than those obtained when \(\lambda\) is assumed to be known, reflecting the cost of pre-estimating the intensity \(\lambda\). 

The deterioration in the convergence rate can be attributed to the increased challenges in controlling the variance of the kernel density estimator. When substituting \(\lambda\) with its estimator, the resulting pair \((X, \hat{\lambda})\) is no longer Markovian nor exponentially ergodic, complicating the analysis. Further insights into these challenges are provided in Remarks \ref{rk: deterioration rate lambda hat} and \ref{rk: pb lambda hat}.

To establish the convergence rates of the proposed estimators, we develop probabilistic tools of broader applicability. Specifically, we employ a Girsanov change of measure to transform the Hawkes process with intensity \(\lambda\) into a Poisson process with constant intensity. To obtain this, we extend the bounds for the exponential moments of the Hawkes process—originally established by \cite{leblanc2024exponential} in the stationary case—to the more general non-stationary setting.

Finally, we validate our theoretical results through numerical simulations, illustrating how the variance of the kernel density estimator depends on the choice of bandwidth parameters.

\paragraph{Outline of the paper.}
The paper is organized as follows: Section \ref{sec:model} introduces the model along with its assumptions. In Section \ref{sec:estim}, we define the kernel density estimator for the invariant density and establish convergence rates, assuming the true invariant distribution belongs to an anisotropic Hölder space and that the intensity \(\lambda\) is known. Section \ref{sec:estLambda} extends this framework by introducing a plug-in kernel estimator for the case where \(\lambda\) is also estimated, and we derive the corresponding convergence rates. Section \ref{s:prob} presents the probabilistic tools underpinning these rates, which may also be of independent interest. In Section \ref{sec:num}, we conduct a simulation study to illustrate the convergence behavior of the variance terms. Finally, Section \ref{sec:proofs} contains the proofs of all main results, while the proofs of the technical results required for the main results are deferred to Section \ref{s: proof technical}.

\section{Model and assumptions}
\label{sec:model}
Let $(\Omega, \mathcal{F}, \mathbb{P})$ be a probability space. We define the Hawkes process for $t \geq 0$ through stochastic intensity representation. We introduce the $1$-dimensional point process 
$N_t$ and its intensity $\lambda_t$, conditional on the history up to time $t$. Precisely, we define $\bar{\mathcal{F}}_t := \sigma (N_s, 0 \leq s \leq t) $ the history of the counting process $N$. The intensity process $(\lambda_t)$ of the counting process $N$ is 
{a} 
$\bar{\mathcal{F}}_t$-predictable process 
that makes {$N_t - \int_0^t \lambda_s ds$}
a $\bar{\mathcal{F}}_t$-local martingale. {We assume that $\lambda$ is left continuous, while $N$ is a cadlag process.}

The jump process $N$ satisfies an initial condition on $(-\infty, 0]$; $N^{\rm in} = N_{(-\infty, 0]}$ is a $\mathcal{F}_0$-measurable boundedly finite point process and the 
 conditional intensity measure of $\left.N\right|_{(0, \infty)}$ is absolutely continuous w.r.t. the Lebesgue measure with density, for $t>0$
\begin{equation*}\label{eq:intensity}
\lambda_t= \xi+\int_{(-\infty, t)} \alpha e^{-\beta (t-u)} {\dd N_u} ,
\end{equation*}
with $\xi>0$ the baseline parameter, $\alpha \in \R^+$ represents the impact of an event on the probability of observing a new jump, and $\beta$ is the decay parameter of the exponential kernel which characterizes the memory of the process. 
We denote the initial condition $\lambda_0$, such that
\begin{equation}{\label{eq: dynamics lambda 1}}
\lambda_t= \xi+   (\lambda_0-\xi)e^{-\beta t}+\int_{0}^{{t-}} \alpha e^{-\beta(t-u)}\dd N_u.
\end{equation}
%
%
The following assumption ensures the nonexplosion of the Hawkes process as $\alpha/\beta$ represents the mean number of children that an immigrant has, where the immigrants arrive according to a Poisson process with intensity $\xi$. 
\begin{ass}\label{ass:hawkes}
$\xi>0$, $ \alpha/\beta < 1$
\end{ass}
In this paper, our goal is to estimate the invariant measure associated with the pair \((X, \lambda)\). To achieve this, it is essential to first confirm the existence of such a measure. For this purpose, we introduce the following assumptions on the coefficients of the stochastic differential equation
\eqref{Eq : model diff X intro}.
\begin{ass}[Assumptions on the coefficients of $X$]\label{ass: X}\text{}
\begin{enumerate}[label=\roman*)]
    \item The coefficients \(b\) and \(\sigma\) are of class \(\mathcal{C}^2\), and there exists a positive constant \(c\) such that, for all \(x \in \mathbb{R}\), we have \( |b'(x)| + |\sigma'(x)| + |a'(x)| \leq c\).
    \item There exist positive constants \(a_1\) and \(\sigma_1\) such that \( |a(x)| < a_1 \) and \( 0 < \sigma^2(x) < \sigma_1^2 \) for all \(x \in \mathbb{R}\).
    \item There exist positive constants \(c'\) and \(q\) such that, for all \(x \in \mathbb{R}\), \( |b''(x)| + |\sigma''(x)| \leq c'(1 + |x|^q)\).
    \item There exist constants \(d \geq 0\) and \(r > 0\) such that, for all \(x\) with \( |x| > r \), we have \( xb(x) \leq -d x^2 \).
\end{enumerate}
\end{ass}
%
%
Following \cite{DLL}, under these assumptions, the pair \((X, \lambda)\) is a Markov process for the general filtration
$\mathcal{F}_t := \sigma(W_s, N_s,0 \le s \le t)$, 
that is ergodic, with an invariant intensity measure \(\pi\). In this paper, we assume that the process is distributed according to this invariant distribution, placing our analysis in the stationary regime:
\[ (X_0, \lambda_0) \sim \pi. \]
{We assume that this measure admits a density with respect to the Lebesgue measure, denoted by \(\pi(\dd x, \dd y) = \pi(x, y)\, \dd x\, \dd y\). In the Poisson case, \cite{bichteler1983clacul} employs Malliavin calculus to characterize conditions on the coefficients ensuring absolute continuity of \(\pi\). Malliavin calculus for Hawkes processes is currently an active research topic \cite{hillairet2021malliavin}; in particular, Theorem 5.9 in \cite{cacitti2025malliavin} provides sufficient conditions for absolute continuity of the law of a Hawkes process. However, this result holds in the absence of a Brownian component, while in our setting both Brownian and Hawkes components are present. To our knowledge, no existing work guarantees the existence of a density in this hybrid case—an interesting open question left for future work. For the purposes of this study, we assume the existence of such a density \(\pi\), which we aim to estimate.}

\section{Estimator of the invariant distribution for a known intensity}
\label{sec:estim}

As mentioned earlier, this section focuses on the estimation procedure in the case where the process \((X_t, N_t)_{t \in [0, T]}\) is continuously observed, under the asymptotic regime where the time horizon \(T\) tends to infinity. We place first the study when the parameter $\theta=(\xi, \alpha, \beta)$ is known. Thus, as $\lambda_t$ is $\mathcal{F}_{t}$-predictable,
its value is known in $t$ as soon as $\theta$ is known together with the jumps of the process $N$, that we know as we observe continuously $(X_t)$.

\subsection{Kernel estimator}

The kernel function is denoted \(K: \mathbb{R} \rightarrow \mathbb{R}\). It satisfies the following properties:
\[
\int_{\mathbb{R}} K(x) \, \dd x = 1, \quad \|K\|_\infty < \infty, \quad \text{supp}(K) \subset [-1, 1], 
\]
We denote \(h_1, h_2 > 0\) as two bandwidth parameters, and define the rescaled kernel functions \(\mathbb{K}_{h_j}\) for \(j = 1, 2\) as:
\[
\mathbb{K}_{h_j}(z) := \frac{1}{h_j} K\left(\frac{z}{h_j}\right), \quad z \in \mathbb{R}.
\]
The kernel density estimator of \(\pi\) at \((x^*, y^*) \in \mathbb{R}^2\) is then given by
\begin{equation}\label{eq:estim}
\widehat{\pi}_{h_1, h_2}(x^*, y^*) :=
\frac{1}{T} \int_{0}^T \mathbb{K}_{h_1}(x^* - X_u)\mathbb{K}_{h_2}(y^* - \lambda_u) \, \dd u.
\end{equation}
The suitability of this estimator for estimating the invariant density is closely tied to the ergodicity of the process. Under appropriate assumptions, \(T^{-1} \int_0^T \delta_{(X_u, \lambda_u)} \, du\) converges to \(\pi\) when \(T \to \infty\). Consequently, \(\hat{\pi}_{h_1, h_2}(x^*, y^*)\) converges almost surely to \(\E[\mathbb{K}_{h_1}(x^* - X_0)\mathbb{K}_{h_2}(y^* - \lambda_0)] \), which converges to \(\pi(x^*, y^*)\) as \(h_1, h_2\) tends to zero.

\subsection{Rate of convergence for Hölder regularity}

We assume that the invariant density $\pi$ belongs to the anisotropic H\"older space \(\mathcal{H}(\beta, \mathcal{L})\), with \(\beta = (\beta_1, \beta_2)\), as defined below, with $ \beta_i > 0.$

\begin{definition}
A function $f$ is in the space $\mathcal{H}(\beta, \mathcal{L})$ with $\beta = (\beta_1, \beta_2)$ and $\mathcal{L} = (\mathcal{L}_1, \mathcal{L}_2)$ if and only if \(f\) is \(p_i\)-times differentiable with respect to \(x_i\), where \(p_i = \left\lfloor \beta_i \right\rfloor\) for \(i \in \{1, 2\}\), and
\[
\left|\frac{\partial^{p_1} f}{\partial x_1^{p_1}}(x_1, x_2) - \frac{\partial^{p_1} f}{\partial x_1^{p_1}}(x_1', x_2)\right| \leq \mathcal{L}_1 |x_1 - x_1'|^{\beta_1 - p_1}, \quad \forall x_1, x_1' \in \mathbb{R}, \forall x_2 \in \mathbb{R},
\]
and
\[
\left|\frac{\partial^{p_2} f}{\partial x_2^{p_2}}(x_1, x_2) - \frac{\partial^{p_2} f}{\partial x_2^{p_2}}(x_1, x_2')\right| \leq \mathcal{L}_2 |x_2 - x_2'|^{\beta_2 - p_2}, \quad \forall x_1 \in \mathbb{R}, \forall x_2, x_2' \in \mathbb{R},
\]
where \(\left\lfloor \beta \right\rfloor\) is the largest integer less than \(\beta\).
\end{definition}

In other words, a function \(f\) belongs to the two-dimensional anisotropic H\"older space defined above if all partial derivatives of \(f\) up to order \(\lfloor \beta_i \rfloor\) are bounded, and the \(\lfloor \beta_i \rfloor\)-th partial derivative is H\"older continuous of order \(\beta_i - \lfloor \beta_i \rfloor\) in the \(i\)-th direction for \(i = 1, 2\).
\begin{ass}\label{ass:reg} Let us assume that 
$\pi\in \mathcal{H}(\beta, \mathcal{L})$, with $\beta=(\beta_1, \beta_2)$. {We suppose that there exists an integer $M \geq \max(\beta_1, \beta_2)$, such that
$$ \forall i \in \{{1}, \dots, M\},~ \int_{\mathbb{R}} K(x) x^i \, \dd x = 0. $$}
\end{ass}
%
It is natural in our context to assume that the invariant density belongs to a H\"older class as described above. 

We aim to analyze this estimator in detail, focusing on the rate of convergence of $\w{\pi}_{h_1,h_2}(x^*, y^*)$ to \(\pi(x^*, y^*)\). To do so, we study the pointwise mean squared error 
and determine its rate of convergence as a function of \(T\). As it is standard, the bound on the mean squared error is derived through the classical bias-variance decomposition,
\begin{equation}{\label{eq: bias-variance}}
 \E[(\hat{\pi}_{h_1, h_2}(x^*, y^*) - \pi(x^*, y^*))^2] = (\E[\hat{\pi}_{h_1, h_2}(x^*, y^*)] - \pi(x^*, y^*))^2 + \text{Var}(\hat{\pi}_{h_1, h_2}(x^*, y^*)).   
\end{equation}
The investigation of these two terms is different. The variance term is more intricate and less standard than the bias.
Let us first study the bias term in the following proposition. 
\begin{prop}\label{prop:biais}
Under Assumptions \ref{ass:hawkes}, \ref{ass: X}, and \ref{ass:reg}, if $\w{\pi}_{h_1, h_2}$ is the estimator given in \eqref{eq:estim} and $y^* \ge \xi$, then there exist constants $C > 0$ such that, for all $T > 0$,
\begin{equation*}
\left(\E[\w{\pi}_{h_1, h_2}(x^*,y^*)]-{\pi}(x^*,y^*)\right)^2 \le C \left( h_1^{2 \beta_1} + h_2^{2 \beta_2}\right).
\end{equation*}
\end{prop}
Note that this bound is identical to those obtained for the bias term in cases where one considers an SDE driven by classical Brownian motion, a Lévy process, or fractional Brownian motion (see Proposition 1 in \cite{comte2013anisotropic}, Proposition 2 in \cite{amorino2021invariant}, and Proposition 5 in \cite{amorino2024fast}, respectively). Given its standard character, we have chosen to defer its proof to Section \ref{s: proof technical}, which is dedicated to the technical results.

We focus now on the variance term. A potentially surprising result is that we can derive two distinct bounds on the variance, depending on whether the pointwise estimation is conducted at the baseline, i.e., whether $y^* = \xi$ or not.
We first present the result for the case where $y^* = \xi$.
\begin{prop}{\label{prop: y=}}
Under Assumptions \ref{ass:hawkes} and \ref{ass: X}, if $\w{\pi}_{h_1, h_2}$ is the estimator given in \eqref{eq:estim} and $y^*=\xi$, then there exist constants $C > 0$ and $T_0 > 0$ such that for $T \ge T_0$,
$$\Var(\w{\pi}_{h_1, h_2}(x^*,\xi)) \le \frac{C}{T} \frac{|\log (h_1 h_2)|}{h_2}.$$
\end{prop}
It is important to note that, throughout, the notation $C$ denotes a general constant, whose value may change from line to line.

Let us now proceed to the variance analysis of the kernel density estimator for the case where $y^* > \xi$. In this scenario, we can establish a sharper bound compared to the previous case, which results in a better convergence rate, see Corollary \ref{coro:rates}.

\begin{prop}{\label{prop: y>}}
Under Assumptions \ref{ass:hawkes} and \ref{ass: X}, if $\w{\pi}_{h_1, h_2}$ is the estimator given in \eqref{eq:estim} and $y^* > \xi$, then for any arbitrarily small $\varepsilon > 0$, there exist constants $C > 0$ and $T_0 > 0$ such that for $T \ge T_0$,
$$\Var(\w{\pi}_{h_1, h_2}(x^*,y^*)) \le \frac{C}{T} \frac{1}{\sqrt{h_2}} \frac{1}{(h_1 h_2)^\varepsilon}.$$
\end{prop}

\begin{rem}{\label{rk: idea proof girsanov}}[On the proofs]
The proofs are relegated in Section \ref{sec:proofs}. For both propositions the proofs rely on examining the covariance of the kernel density estimator at two different time points, $t$ and $s$, and using different bounds depending on the size of $|t - s|$. Specifically, when the size of $|t - s|$ exceeds $\log T$, we leverage the exponential ergodicity of the process $(X_t, \lambda_t)_{t \ge 0}$. For smaller $|t - s|$, the proof diverges between the two cases. When $y^* > \xi$, we employ a change of measure that allows us to transition from our Hawkes process $N$ to a simple homogeneous Poisson process, via a Girsanov change of measure. This approach is also discussed further in Section \ref{s:prob}. In contrast, for $y^* = \xi$, the argument relies on the explicit dynamics of $\lambda_s$, that one can replace in $\mathbb{K}_{h_2}(y^* - \lambda_s) = \mathbb{K}_{h_2}(\xi - \lambda_s)$.
\end{rem}

It is clear from Propositions \ref{prop: y=} and \ref{prop: y>} that we achieve a better convergence rate when the pointwise estimation is conducted at $y^* > \xi$. This appears again, more clearly, when we select the optimal bandwidth $h=(h_1, h_2)$, leading to the convergence rates described in Corollary \ref{coro:rates}. Intuitively, the deterioration in the estimation procedure for $y^* = \xi$ arises because the process never visits the baseline, i.e., the probability of visiting the baseline is zero. Since the process does not reach the baseline, gathering observations and performing pointwise estimation at $y^* = \xi$ becomes challenging.

Mathematically, the method used in the proof of Proposition \ref{prop: y>} is inapplicable to prove Proposition \ref{prop: y=} because it hinges on the condition $y^* - \xi > 0$ (see for example the statement of Lemma \ref{L: g pour q zero}). However, the bounding of the variance for the case $y^* = \xi$ involves a sequence of equalities (see Equation \eqref{eq: y=}), which implies that even though the convergence rate is worse in this case compared to $y^* > \xi$, we believe it is unlikely that the rate from Proposition \ref{prop: y=} could be significantly improved.

Now, using the result of Proposition \ref{prop:biais}, we can easily derive the following convergence rate, based on the bias-variance decomposition given in \eqref{eq: bias-variance}, along with the variance bounds for the cases $y^* = \xi$ and $y^* > \xi$ provided in Propositions \ref{prop: y=} and \ref{prop: y>}, respectively, and the bias bound from Proposition \ref{prop:biais}. The proof of Corollary \ref{coro:rates} also incorporates the rate-optimal choice of bandwidths $h_1$ and $h_2$, and can be found in Section \ref{sec:proofs}.

\begin{coro}\label{coro:rates}   
Suppose Assumptions \ref{ass:hawkes}, \ref{ass: X}, and \ref{ass:reg} hold. Choose $h_2^{\rm opt} = \left(\frac{1}{T}\right)^\frac{1}{2 \beta_2 + 1}$ for $y^* = \xi$ and $h_2^{\rm opt} = \left(\frac{1}{T}\right)^\frac{2}{4 \beta_2 + 1}$ for $y^* >\xi$, with $h^{\rm opt}_1$ small enough to ensure that $h_1^{2 \beta_1}$ is negligible. Then, for any arbitrarily small $\varepsilon > 0$, there exists a constant $T_0 > 0$ such that for all $T \ge T_0$,
\begin{equation*}
\E\left[(\w{\pi}_{h_1^{\rm opt}, h_2^{\rm opt}}(x^*,y^*)-{\pi}(x^*,y^*))^2\right] \lesssim \begin{cases}
\left(\frac{1}{T}\right)^{\frac{2 \beta_2}{2 \beta_2 + 1}} \log(T) \qquad &y^* = \xi, \\
\left(\frac{1}{T}\right)^{\frac{4 \beta_2}{4 \beta_2 + 1} - \varepsilon} \qquad &y^* > \xi.
\end{cases}
\end{equation*}
\end{coro}

Let us comment on the obtained convergence rates by comparison with those established in the literature for the estimation of the invariant density associated with diffusions driven by different types of noise. First, in the classical SDE case, driven by a Brownian motion, it is well known that the convergence rate for estimating the invariant density in the Hölder space \(\mathcal{H}_d (\beta, \mathcal{L})\) is \(T^{- \frac{2\beta}{2\beta + d - 2}}\) for \(d \geq 3\). One can see \cite{dalalyan2007asymptotic, strauch2018adaptive}. This rate has been shown to be optimal in \cite{amorino2024minimax}. The same rate is also obtained in \cite{amorino2021invariant} for a Lévy-driven SDE when \(d \geq 3\).
In lower dimensions, the situation differs: for both classical and Lévy-driven SDEs, the convergence rate does not depend on the regularity. For instance, in the case of \(d = 1\), it is possible to achieve the superoptimal parametric rate \(\frac{1}{T}\), as demonstrated in \cite{castellana1986smoothed} for jump SDE and in \cite{amorino2022optimal} for SDE. For \(d = 2\), the convergence rate becomes \(\frac{\log T}{T}\), as noted in \cite{dalalyan2007asymptotic} and \cite{amorino2021invariant} for SDEs with and without jumps.

In our setting, however, the situation is somewhat different. Although we consider the two-dimensional process and estimate the bidimensional density, we do not achieve the convergence rate \(\frac{\log T}{T}\). The reason for this discrepancy lies in the differing influences of 
$X$ and $\lambda$ on the estimation of the invariant density. 
Specifically, in the dynamics of \(\lambda\), there is no Brownian motion, making it deterministic between jumps. 
This suggests that the convergence rate reflects an intrinsic characteristic of the SDE driven by the Hawkes process in question.

It appears then relevant to compare our findings with those on the estimation of the invariant density of a stochastic two-dimensional hypo-elliptic system, as studied in \cite{delattre2022rate}. Indeed, both results share notable similarities, including a dichotomy in the pointwise rate of convergence depending on the point \(y^*\) where the estimation procedure is performed. Furthermore, in \cite{delattre2022rate}, the final convergence rate depends solely on \(\beta_1\) or \(\beta_2\), depending on the relative positions of these two parameters (see Theorems 1 and 2 therein), which aligns with our findings. These parallels arise because our system also exhibits features of hypoelliptic diffusion, as the stochastic differential equation governing \(\lambda\) is degenerate, lacking a Brownian component.

\subsection{Comments}{\label{s: comments}}

\paragraph{On the minimax optimality.}
One may wonder whether it is possible to improve the convergence rate we have derived above by proposing an alternative density estimator. Specifically, one might ask if the estimator we propose is optimal in a minimax sense. While we believe the variance bounds we provide are sharp and lead to convergence rates that are unlikely to be further improved, we are currently unable to definitively address the question of optimality, as we have yet to establish a lower bound for the mean squared error. It is worth noting that in frameworks analogous to ours, the proof of lower bounds is often based on the two-hypothesis method, which relies heavily on the link between the drift coefficient and the invariant density (see, for example, \cite{amorino2021rate}, \cite{amorino2024minimax} and \cite{delattre2022rate}, where lower bounds are proven for diffusions without and with jumps, and for hypoelliptic diffusions, respectively). This link is straightforward when the invariant density can be explicitly expressed, which is not the case here. Although it may be possible to establish such a connection using the adjoint of the generator, as in Section 6.1 of \cite{amorino2021rate}, it is important to note that we are now dealing with the generator associated with the two-dimensional diffusion \((X_t, \lambda_t)\). This problem is significantly more challenging than the one addressed in \cite{amorino2021rate} and requires more space than we can dedicate in this paper. Thus, we leave this as a topic for future research.

\paragraph{On the independence from $\beta_1$.} 
It is particularly surprising that our convergence rates are entirely independent of \(\beta_1\). This occurs because we can separate the contributions of the jump intensity from that of the diffusion. In our case, the jump-diffusion \(X\) is one-dimensional, allowing us to achieve a superoptimal convergence rate 
which remains independent of the regularity \(\beta_1\). 
However, this usually changes when considering a higher-dimensional diffusion \(X\), where dependence on \(h_1\) (and thus \(\beta_1\)) appears. We believe this dependency would also arise in our framework if extended to higher dimensions, explaining why our current convergence rate depends only on \(\beta_2\). We have chosen to focus here on the one-dimensional case for \(X\), the higher-dimensional scenario is left for future investigation.

\paragraph{On possible extensions.}
We conclude this section by discussing some other potential extensions of the results stated above. In particular, our estimators are based on the continuous observation of the paths of the processes \((X, \lambda)\). However, we claim that the results in this section still hold if only discrete observations of the process are available, provided the estimator in \eqref{eq:estim} is adapted to the discrete framework, as done in \cite{amorino2023estimation} or \cite{maillet2024estimation}. The impact of discretization becomes relevant, however, when the jump intensity \(\lambda\) is unknown and must be pre-estimated, as discussed in Section \ref{sec:estLambda}.

Moreover, it is evident from our main results that the optimal bandwidth \(h_2\) depends on the unknown smoothness parameter \(\beta_2\). Therefore, it may be advantageous to propose an adaptive procedure, similar to the one introduced by Goldenshluger and Lepski in \cite{goldenshluger2011bandwidth}, which allows bandwidth selection based solely on the data, without requiring prior knowledge of \(\beta_2\). While this approach has been studied for the estimation of invariant densities in other models, the analogous procedure in our framework remains unexplored and is left for future investigation.

\subsection{Estimation of the invariant density of the intensity}{\label{s: estim int lambda}}
This section focuses on estimating the invariant density of the jump intensity, \(\lambda\), which can be derived as a by-product of the results from the previous section. This is particularly noteworthy as, up to a logarithmic factor, it demonstrates the possibility of achieving a superoptimal parametric rate, see Corollary \ref{c: rate lambda} and following comments.

Let us consider a natural modification of the kernel density estimator proposed in \eqref{eq:estim} for estimating the invariant density of the jump intensity \(\lambda_t\) alone. For \(y^* \in \mathbb{R}\), this estimator is given by 
\begin{equation}{\label{eq: estim lambda}}
\widehat{\pi}_{h_2}(y^*) :=
\frac{1}{T} \int_{0}^T \mathbb{K}_{h_2}(y^* - \lambda_u) \, \mathrm{d} u.    
\end{equation}
Similar to the previous estimator $\w{\pi}_{h_1, h_2}$,
two different convergence rates can be achieved depending on whether \(y^* = \xi\). 

What is particularly striking is that, in the best-case scenario where \(y^* > \xi\), it is possible to obtain a variance bound that is independent of the bandwidth, resulting in a superoptimal convergence rate.

\begin{prop}\label{prop: var lambda}
Under Assumptions \ref{ass:hawkes} and \ref{ass: X}, if \(\hat{\pi}_{h_2}\) is the estimator given in \eqref{eq: estim lambda}, then for any arbitrarily small $\varepsilon > 0$, there exist constants \(C > 0\) and \(T_0 > 0\) such that for \(T \ge T_0\), 
\[
\Var(\hat{\pi}_{h_2}(y^*)) \le    
\begin{cases}
\frac{C}{T} \frac{|\log h_2|}{h_2} & \text{if } y^* = \xi, \\
\frac{C}{T}(|\log h_2| + h_2^{- \varepsilon}) & \text{if } y^* > \xi.
\end{cases}
\]
\end{prop}

The proof of Proposition \ref{prop: var lambda} is provided in Section \ref{sec:proofs} and follows a similar approach to the bounds established in Propositions \ref{prop: y=} and \ref{prop: y>}. However, the absence of the diffusion term significantly simplifies the analysis. In particular, it allows us to consider a different, bounded function \(f_1\) in the change of measure discussed in Section \ref{s:prob}, resulting in a tighter bound. 

The final convergence rate for the estimation of the invariant density of \(\lambda\) is summarized in the following corollary. Let us note that 
the bias term is controlled as in Proposition \ref{prop:biais}.
%
\begin{coro}\label{c: rate lambda}   
Suppose Assumptions \ref{ass:hawkes}, \ref{ass: X}, and \ref{ass:reg} hold. Then, under the rate-optimal choice of bandwidths, for any arbitrarily small \(\varepsilon > 0\), there exist constants \(C > 0\) and \(T_0 > 0\) such that for all \(T \ge T_0\), for $h_2^{\rm opt}= T^{-1/(2\beta_2+1)}$ for $y^*=\xi$ and $h_2^{\rm opt}=  T^{-1/(2\beta_2+\varepsilon)}$ for $y^*>\xi$ 
\[
\E[(\hat{\pi}_{h_2^{\rm opt}}(y^*) - \pi(y^*))^2] \le 
\begin{cases}
C \left(\frac{1}{T}\right)^{\frac{2 \beta_2}{2 \beta_2 + 1}} \log(T) & \text{if } y^* = \xi, \\
\frac{C}{T}(\log T + T^{\varepsilon}) & \text{if } y^* > \xi.
\end{cases}
\]
\end{coro}

A few remarks on these convergence rates are necessary. In \cite{krell2021nonparametric}, the authors study nonparametric estimation of jump rates for a specific class of one-dimensional piecewise deterministic Markov processes (PDMP), constructing an estimator for the stationary density. Up to a logarithmic factor, their convergence rate matches what we found for the case \(y^* = \xi\) (see Proposition 6 of \cite{krell2021nonparametric} and their subsequent comments). 
It is important to highlight that, for our specific PDMP process \(\lambda\), this represents only the worst-case scenario. In the more favorable case \(y^* > \xi\), we can estimate the jump intensity at a superoptimal rate. It is natural to wonder whether our findings could be extended to more general PDMPs, this is so far an open question, but one that certainly deserves further exploration.



\section{Estimator of the invariant distribution for an unknown intensity}

\label{sec:estLambda}

In the previous section, we introduced an estimator based on the continuous observation of the entire path of the pair $(X_t, \lambda_t)$, for $t \in [0, T]$. However, the assumption that both the process and the jump intensity $\lambda_t$ can be observed may seem unrealistic. To address this, we now weaken this assumption, considering only the observation of the process $X$. Since observing $X$ continuously also reveals the jumps, the parameter $\theta = (\xi, \alpha, \beta)$ can still be estimated. Based on these estimates, we can then construct an estimator for the jump intensity $\lambda$.

\subsection{Plug-in kernel estimator}\label{sec:estimTheta}
Let $\w{\theta}_{{T}}$ denote the maximum likelihood estimator proposed in \cite{Ogata} and \cite{clinet2017statistical}. In these works, the authors estimate the parameters $\theta$ using the maximum likelihood estimator $\w{\theta}_{{T}}= (\w{\xi}, \w{\alpha}, \w{\beta})$. Thanks to Theorem 3.14 in \cite{clinet2017statistical}, we know that for any $k \geq 1$,
\begin{equation}{\label{eq: parametric rate 1.5}}
  \E\left[ \big|\w{\theta}_{{T}}- \theta\big|^k\right] \leq \frac{C}{{T}^{k/2}}.  
\end{equation}
Let us assume we dispose of a second dataset, independent of the first, that we use in order to estimate the parameter $\theta$ as above. Then, by substituting the estimated parameters into the dynamics of $\lambda$, we obtain an estimator for the jump intensity process $\lambda_t$ for any $t \in [0, T]$. Specifically, we define the estimated intensity process as follows:
\begin{equation}{\label{eq: def lambda hat 2}}
  \w{\lambda}_t=\w{\xi}+ \sum_{T_i < t} \w{\alpha} \exp(-\w{\beta}(t-T_i) )+ (\wt{\lambda}_0-\w{\xi})e^{-\w{\beta} t}, 
\end{equation}
where $\wt{\lambda}_0$ is an approximation of the initial condition $\lambda_0$.
We can then show that the parametric rate from \eqref{eq: parametric rate 1.5} for the estimation of $\theta$ extends to the estimation of the jump intensity process as well, provided that the initial condition is approximated properly. This is formalized in the following theorem, whose proof can be found in Section \ref{sec:proofs}.

\begin{theo}{\label{th: bound lambda hat}}
Consider $\lambda_t$ and $\widehat{\lambda}_t$ as defined in \eqref{eq: dynamics lambda 1} and \eqref{eq: def lambda hat 2}, respectively. Assume there exists $\varepsilon \in (0,1)$ such that $\beta > \varepsilon$ and $\widehat{\beta} > \varepsilon$. Then, for any $q \ge 2$ and any $t \in [0,T]$, there exists a constant $C > 0$ such that
\begin{equation}\label{eq:boundlambda}
\E[|\widehat{\lambda}_t - \lambda_t|^q]^\frac{1}{q} \le \frac{C}{\sqrt{{T}}} + \E[|\lambda_0 - \widetilde{\lambda}_0|^q]^\frac{1}{q} e^{- \varepsilon t}.
\end{equation}
\end{theo}
It is important to note that Theorem \ref{th: bound lambda hat} 
holds intrinsic interest. Specifically, it shows that the estimator $\widehat{\lambda}_t$ proposed in \eqref{eq: def lambda hat 2} performs well, achieving a parametric rate of estimation for $\lambda_t$ for any $t \in [0, T]$—provided the initial condition is approximated sufficiently well so that the second term in the bound above is negligible compared to the first. 
As we do not have full control over the quality of the approximation $\wt{\lambda}_0$ for $\lambda_0$, the bound motivates us to consider $t$ large enough to ensure that the second term is negligible. Specifically, we want $e^{-\varepsilon t} \le \frac{C}{\sqrt{T}}$, which leads to the condition:
\[
t \ge \Tm
\]
with
\begin{equation}\label{eq:Tm}
  \Tm := \frac{1}{2 \varepsilon} \log(T).
\end{equation}

With the aim of proposing a kernel density estimator for the invariant density based on the observation of $X$, we begin by using data only from time $\Tm$ onward. Notably, a natural estimator for $\pi(x^*,y^*)$ when only $X$ is observed consists of replacing $\lambda_t$ with $\widehat{\lambda}_t$. This leads to the following estimator:
\begin{equation}\label{eq:estim2}
\widetilde{\pi}_{h_1, h_2}(x^*, y^*) := \frac{1}{T-\Tm} \int_{\Tm}^T \mathbb{K}_{h_1}(x^* - X_u)\mathbb{K}_{h_2}(y^* - \widehat{\lambda}_u) \, du.
\end{equation}

\subsection{Rates of convergence for the plug-in estimator}
We now aim to analyze the performance of this new estimator by studying the convergence rate of its pointwise mean squared error.
Specifically, we seek to bound $\E[(\wt{\pi}_{h_1,h_2}(x^*,y^*) - \pi(x^*,y^*))^2]$. 

There are two possible approaches for achieving this. The first involves bounding the error incurred when transitioning from $\wt{\pi}_{h_1,h_2}$ to $\widehat{\pi}_{h_1, h_2}$ and then applying the results from Section \ref{sec:estim}.
The second approach directly bounds the mean squared error of $\wt{\pi}_{h_1,h_2}$ using a bias-variance decomposition.
We opt for the first approach, as the second presents several challenges due to the fact that the pair $(X_t, \widehat{\lambda}_t)$ is neither Markov nor satisfies the exponential ergodicity property, see Remark \ref{rk: pb lambda hat} for details.

We can assume that the estimator \(\widehat{\pi}_{h_1, h_2}\) has been constructed solely from the observations of \((X_t, \lambda_t)\) for \(t \in [\Tm, T]\), as this will simplify our computations.

We decompose the pointwise mean squared error as follows:
\begin{eqnarray}
\label{eq:decompHatlambda}
\E\left[(\wt{\pi}_{h_1,h_2}(x^*,y^*) - \pi(x^*,y^*))^2\right] &\leq& 2 \E\left[(\wt{\pi}_{h_1,h_2}(x^*,y^*) - \widehat{\pi}_{h_1,h_2}(x^*,y^*))^2\right] \nonumber\\
&&+ 2\E\left[(\widehat{\pi}_{h_1,h_2}(x^*,y^*) - \pi(x^*,y^*))^2\right].
\end{eqnarray}
Corollary \ref{coro:rates} provides control over the second term, leaving us to bound the first. Naturally, this control is achieved by leveraging the bound on the error incurred when approximating $\lambda_t$, as stated in Theorem \ref{th: bound lambda hat}. Interestingly—and somewhat unsatisfactorily—despite $\widehat{\lambda}_t$ approximating $\lambda_t$ at a parametric rate, the error in moving from $\wt{\pi}_{h_1,h_2}$ to $\widehat{\pi}_{h_1,h_2}$ is not negligible. This results in a deterioration of the convergence rate compared to the one derived in Corollary \ref{coro:rates}, see Remark \ref{rk: deterioration rate lambda hat} for details. This deterioration is the cost of not knowing $\lambda$, as made clear in the following proposition.

\begin{prop}{\label{prop: Ta}}
Suppose Assumptions \ref{ass:hawkes}, \ref{ass: X}, and \ref{ass:reg} hold. Assume there exists $\tilde{\varepsilon} \in (0,1)$ such that $\beta > \tilde{\varepsilon}$ and $\widehat{\beta} > \tilde{\varepsilon}$. Then, for any arbitrarily small $\varepsilon > 0$, there exist constants $C > 0$ and $T_0 > 0$ such that, for all $T \ge T_0$:
\begin{equation}{\label{eq: error lambdahat}}
\E[(\w{\pi}_{h_1, h_2}(x^*,y^*) - \wt{\pi}_{h_1, h_2}(x^*,y^*))^2]  \le C \frac{1}{h_1^{1 + \varepsilon} h_2^4} \frac{1}{T}.    
\end{equation}
Thus, the rate-optimal choice for $h_1^{\rm opt}$ and $h_2^{\rm opt}$,  leads to:
\begin{equation}{\label{eq: rate lambdahat}}
 \E[({\wt{\pi}}_{h_1^{\rm opt}, h_2^{\rm opt}}(x^*,y^*) - {\pi}(x^*,y^*))^2]  \le C \left(\frac{1}{T}\right)^{\frac{2 \beta_2}{\beta_2(2 + \frac{1}{\beta_1}) + 4} - \varepsilon}.   
\end{equation}
\end{prop}
The proof of this result can be found in Section \ref{sec:proofs}, the choices $(h_1^{\rm opt}, h_2^{\rm opt})$ are given in Equation \eqref{eq: h lambdahat}.
The first observation to make regarding the rate above is that, unlike the results in Section \ref{sec:estim}, the convergence rate in \eqref{eq: rate lambdahat} no longer depends on whether \( y^* = \xi \). It now also depends on \(\beta_1\), a consequence of \eqref{eq: error lambdahat} involving \(h_1\). Even when the regularity \(\beta_1\) is large, the rate in \eqref{eq: rate lambdahat} is worse compared to both rates presented in Corollary \ref{coro:rates}. This becomes even clearer when comparing the bound in \eqref{eq: error lambdahat} with those on the variance of \(\w{\pi}\), which makes the variance terms negligible in comparison. Consequently, the trade-off between \(h_1\) and \(h_2\) is determined by balancing the condition in \eqref{eq: error lambdahat} with the bias bound from Proposition \ref{prop:biais}, which appears in the second term of \eqref{eq:decompHatlambda}.

We believe that this deterioration in the convergence rate reflects the cost of extending the results to the case where the jump intensity is unknown. While we do not claim that the rate found in \eqref{eq: rate lambdahat} is optimal, this framework introduces many new challenges, see Remarks \ref{rk: deterioration rate lambda hat} and \ref{rk: pb lambda hat} for further details. Whether this bound can be improved remains an open question for future research.

\begin{rem}{\label{rk: deterioration rate lambda hat}}
One might wonder why the rate in \eqref{eq: rate lambdahat} is unsatisfactory and not negligible compared to the bounds from the previous section, despite starting with a strong approximation of \(\lambda_t\).
The reason lies in the nature of our kernel density estimator. Indeed, when transitioning from the kernel applied to \(\lambda\) to one applied to \(\widehat{\lambda}\), an additional factor of \(h_2\) emerges due to the need to differentiate \(\mathbb{K}_{h_2}\). This factor alone, however, does not fully explain the deterioration of the rate, nor the new dependence on \(h_1\) and \(\beta_1\).

The root cause of this issue stems from the fact that the technique used in the previous section to separate the contribution of \(X\) no longer applies here. In the earlier case, it was possible to isolate the effect of \(X\) in the kernel by conditioning the expectation on the path of \(\lambda\). In the current setting, attempting the same does not succeed in isolating the contribution of the kernel of \(X\) because \(\w{\lambda}_t\) is not measurable with respect to the \(\sigma\)-algebra generated by \(\lambda\).   
\end{rem}

\begin{rem}{\label{rk: pb lambda hat}}
Motivated by the deterioration of the rate discussed above,
one might suspect that the convergence rate obtained in \eqref{eq: rate lambdahat} could be improved by directly applying a bias-variance decomposition to the estimator \(\w{\pi}_{h_1, h_2}\), similar to Section \ref{sec:estim}. However, it is crucial to note that this methodology would no longer apply here. We believe that even though the bias bound might not hold anymore, it could potentially be replaced by a comparable result with negligible terms. This relies on knowing that the maximum likelihood estimator satisfies \(|\E[ \w{\theta}_{{T}} - \theta]| \leq \frac{C}{{T}^{\gamma}}\) for some \(\gamma > \frac{1}{2}\), which seems like a reasonable assumption. 
The primary challenge, however, lies in analyzing the variance term. Since the pair \((X, \widehat{\lambda})\) is no longer Markovian nor exponentially ergodic, this poses a significant issue. As explained in Remark \ref{rk: idea proof girsanov}, the variance bound heavily depends on the ergodicity of the pair, making it problematic here. One might speculate that exponential ergodicity could still be established for the pair \((X, \widehat{\lambda})\), since \(\widehat{\lambda}\) is constructed from the same Hawkes process \(N\) as \(\lambda\), allowing the \(\sigma\)-algebra generated by \((X, \lambda)\) to be replaced by that generated by \((X, N)\). While this might seem to address the issue on one front, the fact that \(\widehat{\lambda}\) is not measurable with respect to this new \(\sigma\)-algebra—due to the additional randomness introduced by the estimator \(\widehat{\theta}_T\)—renders this approach unworkable.
\end{rem}


\section{Probabilistic results}
\label{s:prob}
This section is dedicated to presenting the probabilistic results required to establish our main findings. We chose to devote a separate section to these results, as they may hold independent interest. Specifically, we need to bound the conditional expectation with respect to \(\lambda_0\) of a certain functional $f$ being $\mathcal{F}_s$-measurable for some possibly large $s$. A key aspect of this is understanding how many jumps the process has experienced before the time of interest $s$. 
Suppose we are given some time \(s\), and the process has experienced \(q\) jumps before \(s\) (i.e., \(N_s = q\)). We then examine in detail the jump times at which these jumps occurred, denoted as \(T_1, \ldots, T_q\). It is important to note that these jump times are random variables and are not independent of one another. The function 
can thus be viewed as a function of the time increments \(s - T_q\), \(T_q - T_{1}\), \(\ldots\), \(T_2 - T_1\), and \(T_1\). We are particularly interested in the case where \( f \) can be factorized multiplicatively, meaning that \( f \) on 
the event $\{N_s = q\}$
 can be expressed as \( f^{(q)}(s - T_q, T_q - T_1, \ldots, T_2 - T_1, T_1) = f_1(s - T_q) f_2^{(q)}(T_q - T_1, \ldots, T_2 - T_1, T_1) \). Here, \( f_1 \) is integrable near 0, while \( f_2^{(q)} \) is bounded and has compact support (see the theorem below for further details). 
 With these notations $f$ is given by 
 \begin{align*}
 f&=\sum_{q=1}^\infty f^{(q)}(s - T_q, T_q - T_1, \ldots, T_2 - T_1, T_1) \one_{N_s = q}
 \\&=
 \sum_{q=1}^\infty f_1(s - T_q) f_2^{(q)}(T_q - T_1, \ldots, T_2 - T_1, T_1) \one_{N_s = q}.
 \end{align*}
  Note that, in our application, \( f_2^{(q)} \) will be an indicator function.



We aim to establish bounds on the conditional expectation of this quantity, as described in Theorem \ref{th: girsanov} below.

\begin{theo}{\label{th: girsanov}}
Let \(N\) be a self-exciting Hawkes process with intensity \((\lambda_t)_{t \ge 0}\), following the dynamics given in \eqref{eq: dynamics lambda 1}. Consider a function \(f\) as defined above, with \(f_1\) a real positive decreasing function, integrable in the neighborhood of $0$ and $f_2^{(q)}: (\mathbb{R}^+)^q \rightarrow \mathbb{R}_ +$ bounded and such that $f_2^{(q)}(t_q, \ldots, t_1) \neq 0$ only if \(t_1 \in J^{(q)} = J^{(q)}(t_2, \ldots, t_q)\), where \(J^{(q)}\) is some interval in \(\mathbb{R}^+\). We assume that $\sup_{q \ge 1} |J^{(q)}| \le \overline{J} <\infty$.

Then, for all $\varepsilon > 0$, there exists a constant \(C_\varepsilon > 0\) such that, for all $s \ge 0$
\[
\sum_{q=1}^{\infty}A^{(q)} := \sum_{q=1}^{\infty} \mathbb{E}[f^{(q)}(s - T_q, \ldots, T_2 - T_1, T_1) \one_{N_s = q} \mid \lambda_0 = y_0] \leq C_\varepsilon \left(\int_0^{\overline{J}} f_1(u) \, du \right)e^{\varepsilon s},
\]  
where the constant $C_\varepsilon$ is uniform for $y_0$ in a compact interval of $\R$.
\end{theo}

The rigorous proof of Theorem \ref{th: girsanov} is provided in Section \ref{s: proof technical}. Here, we offer a heuristic outline of the proof along with the statements of key probabilistic results necessary to derive it, as we believe these may be of independent interest.

The main idea is that if the process \(N\) were simply a homogeneous Poisson process, it would be possible to obtain a sharp bound for \(A^{(q)}\) by exploiting the independence between \(T_1\) and \(T_j - T_1\), and by integrating with respect to the density of \(T_1\). To facilitate this, we introduce a Girsanov change of measure, which transforms the Hawkes process \(N\) with intensity \(\lambda\) into a Poisson process with constant intensity \(\xi\). 
For this purpose, we define, for all \(s \geq 0\), the random variable
\begin{equation} \label{eq: def L}
	\mathbb{L}_s =\left( \prod_{0 \le T_j \le s} \frac{\lambda_{T_j}}{\xi} \right) \times \exp\left(- \int_0^s (\lambda_u - \xi) \, du \right).
\end{equation}
If we define \(\mathbb{Q}\) as the probability measure such that 
\[
\frac{d\mathbb{Q}}{d\mathbb{P}} = \frac{1}{\mathbb{L}_s} \quad \text{on } \mathcal{F}_s,
\]
then the process \((N_u)_{0 \leq u \leq s}\) becomes a Poisson process with intensity \(\xi\) under \(\mathbb{Q}\) (see \cite{jacodMultivariatePointProcesses1975a}).

Observe that, thanks to this change of variables, we can express 
\[
A^{(q)} = \mathbb{E}_{\mathbb{Q}}[f^{(q)}(s - T_q, \ldots, T_2 - T_1, T_1) \one_{N_s = q} \mathbb{L}_s \mid \lambda_0 = y_0].
\]
Our goal is to eliminate the dependence on \(T_1\). {As we are going to use a shift of time, it will be convenient to work with $\mathbb{L}_{T_1 + s}$ rather than with $\mathbb{L}_s$. We therefore replace this quantity in $A^{(q)}$, remarking we are allowed as $\mathcal{F}_s \subset \mathcal{F}_{T_1 + s}$.} To eliminate such dependence, we note that \(\lambda_{T_1} \in [\xi , \lambda_0 ]\) and introduce upper and lower bounds for \(\lambda_u\), denoted \(\overline{\lambda}_u\) and \(\underline{\lambda}_u\), respectively, such that \(\underline{\lambda}_u \leq \lambda_u \leq \overline{\lambda}_u\) for \(u > T_1\), and \(\overline{\lambda}_u = \underline{\lambda}_u = \lambda_u\) for \(u \leq T_1\), see the proof of Theorem \ref{th: girsanov}, particularly \eqref{eq : lambda bar sur} and \eqref{eq : lambda bar sous}, for more details. As a result, we can define an upper bound for the process \(\mathbb{L}_s\), eliminating the dependence on \(\lambda_{T_1}\) by setting
\begin{equation} \label{eq: def overline L}
	\overline{\mathbb{L}}_s 
	=\left( \prod_{0 \leq T_j \leq s} \frac{\overline{\lambda}_{T_j}}{\xi} \right) \times
	\exp\left(- \int_0^s (\underline{\lambda}_u - \xi) du \right).
\end{equation}

To further eliminate the dependence of the first $q-1$ coordinates of \(f_2^{(q)}\) on \(T_1\), we define the process \(\widetilde{N}_u = N_{T_1 + u} - N_{T_1}\), resetting the counting process to 0 after \(T_1\). Note that under \(\mathbb{Q}\), \(N\) is a Poisson process, and thus \(\widetilde{N}_u\) is also a Poisson process, independent of \(\mathcal{F}_{T_1}\). We denote the jump times of \(\widetilde{N}\) as \(\widetilde{T_j} = T_{j+1} - T_1\) for \(j \geq 0\). We can now express \(A^{(q)}\), using \(\mathbb{L}_{T_1 + s}\), as a change of measure for the shifted process \(\widetilde{N}\):
\[
A^{(q)} \leq \mathbb{E}_{\mathbb{Q}}[f_1(s - \widetilde{T}_{q-1} - T_1) f_2^{(q)}( \widetilde{T}_{q-1}, \ldots, \widetilde{T_1}, T_1) \one_{N_s = q} \overline{\mathbb{L}}_{T_1 + s} \mid \lambda_0 = y_0].
\]

Next, we introduce the \(\sigma\)-algebra generated by the shifted Hawkes process, denoted \(\widetilde{\mathcal{F}}_u := \sigma(\widetilde{N}_v, v \leq u)\). In the proof, we will rigorously show {(see \eqref{eq: 54.5} and the proof there above)} that there exists a constant \(C\), depending on \(\xi\), \(\lambda_0\), and \(\alpha\), such that
\begin{equation}{\label{eq:bound bar tilde L}}
{{\mathbb{L}}_{T_1 +s}} \leq \overline{\mathbb{L}}_{T_1 + s} \leq \widetilde{\mathbb{L}}_s \frac{\lambda_0}{\xi},   
\end{equation}
where \(\widetilde{\mathbb{L}}_s\) is \(\widetilde{\mathcal{F}}_s\)-measurable.

We can now write \(\{N_s = q\}\) as \(\{T_1 \leq s - \widetilde{T}_{q-1}\} \cap \{\widetilde{N}_s = q - 1\}\), and use \(\widetilde{T}_j = T_{j+1} - T_1\) to deduce
\begin{equation}{\label{eq: meas tilde F}}
A^{(q)} \leq \mathbb{E}_{\mathbb{Q}}\left[\mathbb{E}_{\mathbb{Q}}\left[f_1(s - \widetilde{T}_{q-1} - T_1) f_2^{(q)}( \widetilde{T}_{q-1}, \ldots, \widetilde{T_1}, T_1) \one_{T_1 \leq s - \widetilde{T}_{q-1}} \mid \widetilde{\mathcal{F}}_s \right] \one_{\{\widetilde{N}_s = q-1\}} \widetilde{\mathbb{L}}_s \mid \lambda_0 = y_0\right] \frac{y_0}{\xi}.
\end{equation}

Since \(\widetilde{\mathbb{L}}_s\) and \(\widetilde{N}_s\) are \(\widetilde{\mathcal{F}}_s\)-measurable, the inner conditional expectation depends on \((\widetilde{T}_j)_{j = 1, \dots, q-1}\), which are also \(\widetilde{\mathcal{F}}_s\)-measurable due to the restriction \(\widetilde{T}_{q-1} \leq s\). Consequently, this conditional expectation reduces to an integral over the conditional law of \(T_1\) given \(\widetilde{\mathcal{F}}_s\). Under \(\mathbb{Q}\), the process \(N\) is a Poisson process, making \(T_1\) independent of \(\widetilde{\mathcal{F}}_s = \sigma(N_{T_1 + u} - N_{T_1}, 0 \leq u \leq s)\), and the conditional law of \(T_1\) is exponential with parameter \(\xi\).

Thus, we have the bound
\begin{multline*}
\mathbb{E}_{\mathbb{Q}}\left[f^{(q)}(s - \widetilde{T}_{q-1} - T_1, \widetilde{T}_{q-1}, \ldots, \widetilde{T_1}, T_1) \one_{T_1 \leq s - \widetilde{T}_{q-1}} \mid \widetilde{\mathcal{F}}_s \right] \\
\leq \int_0^{s - \widetilde{T}_{q-1}} f_1(s - \widetilde{T}_{q-1} - u) f_2^{(q)}(\widetilde{T}_{q-1}, \ldots, \widetilde{T_1}, u) \xi e^{- \xi u} du.
\end{multline*}
Under the theorem's hypothesis on \(f_2^{(q)}\) and using that $f_1$ is decreasing, this is bounded by \(c \int_0^{\overline{J}} \xi f_1(v) \dd v\). Using this bound, along with {\eqref{eq: meas tilde F}}, we conclude
\begin{equation}{\label{eq: end prob}}
\sum_{q = 1}^\infty A^{(q)} 
\leq 
c \left(\int_0^{\overline{J}} f_1(v) \dd v\right) \mathbb{E}_{\mathbb{Q}}[\overline{\mathbb{L}}_{T_1 + s} \mid \lambda_0 = y_0].    
\end{equation}
The proof is completed by bounding the conditional expectation above. Recall that, since \(\mathbb{L}_{T_1 + s}\) is a Radon-Nikodym derivative, \(\mathbb{E}_{\mathbb{Q}}[\mathbb{L}_{T_1 + s}] = \mathbb{E}_{\mathbb{P}}[1] = 1\). The remaining task is to analyze the deviations between \(\mathbb{L}_{T_1 + s}\) and \(\overline{\mathbb{L}}_{T_1 + s}\), which is done through Lemma \ref{L: upper bound Girsanov L1} below, with the full proof presented in Section \ref{s: proof technical}.

\begin{lemme} \label{L: upper bound Girsanov L1}
	For all $\varepsilon>0$, there exists $C_\varepsilon>0$, such that, for all $s \ge 0$,
\begin{equation*}
	\mathbb{E}_{\mathbb{Q}} \left[\overline{\mathbb{L}}_{T_1+s} \mid \lambda_0=y_0\right] \le 
	C_\varepsilon e^{\varepsilon s}
	\end{equation*}
	where the constant $C_\varepsilon$ is uniform for $y_0$ in a compact of $\R$.
\end{lemme}

{Remark that, in the statement above, one is not allowed to let $\epsilon$ go to zero, as in such case the constant $C_\epsilon$ explodes. However, for any fixed $\epsilon > 0$, the statement holds, and this is enough to obtain our main results.}

The proof is moreover based on a bound for the exponential moments of the Hawkes process. The stationary case is studied in Theorem 1 of \cite{leblanc2024exponential}. We extend it to the non-stationary case in Lemma \ref{L: extension Leblanc}. 
To extend the result to the non-stationary version, we construct a coupling between two Hawkes processes by leveraging the cluster representation of the process, as described in \cite{hawkesClusterProcessRepresentation1974} and \cite{mollerPerfectSimulationHawkes2005c}.
\begin{lemme}
 \label{L: extension Leblanc}
Let $\tilde{N}$ be an exponential Hawkes process with intensity $\widetilde{\lambda}$ solution of \eqref{eq: dynamics lambda 1} satisfying Assumption \ref{ass:hawkes}.  Then, for any $\displaystyle K \in \left(1,\left(\frac{\beta}{2\alpha}+\frac{1}{2}\right)^{1+\frac{2}{1-\alpha/\beta}} \right]$ and $t \ge 0$
\begin{equation}\label{eq : majo moment expo Hawkes}
\E\left[K^{\widetilde{N}([0,t])}\mid \widetilde{\lambda}_0=y_0\right] \le \exp \left(
		\big(K^{1+\frac{2}{1-\beta/\alpha}} -1\big)\times t \times (\min(y_0,\xi))
		\right).
\end{equation}    
\end{lemme}	
The proof can be found in Section \ref{s: proof technical}.

\section{Numerical illustration}
\label{sec:num}

\begin{figure}[hbtp]
    \centering
\includegraphics[width=0.4\linewidth]{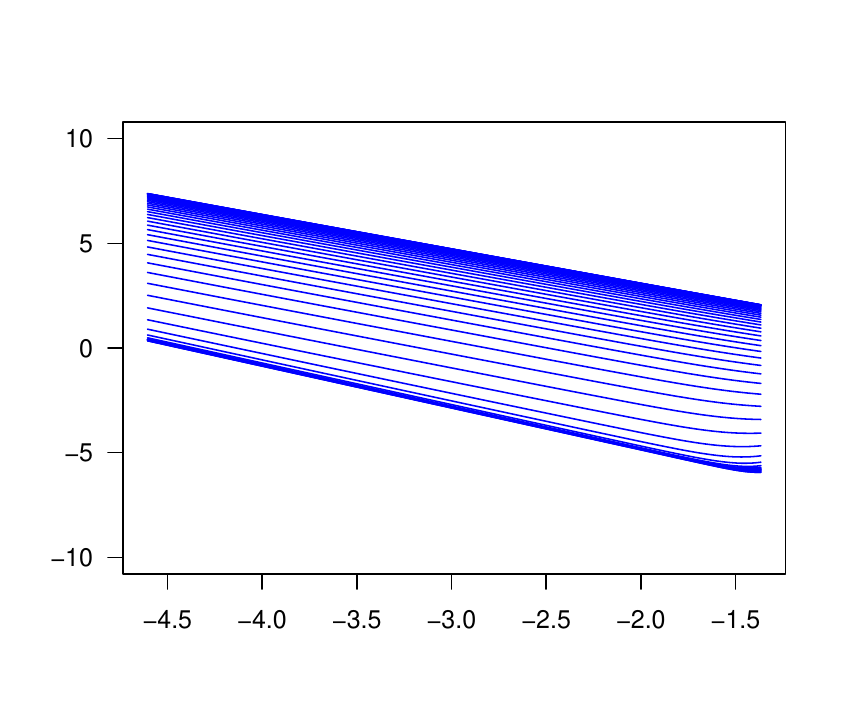}
  \includegraphics[width=0.4\linewidth]{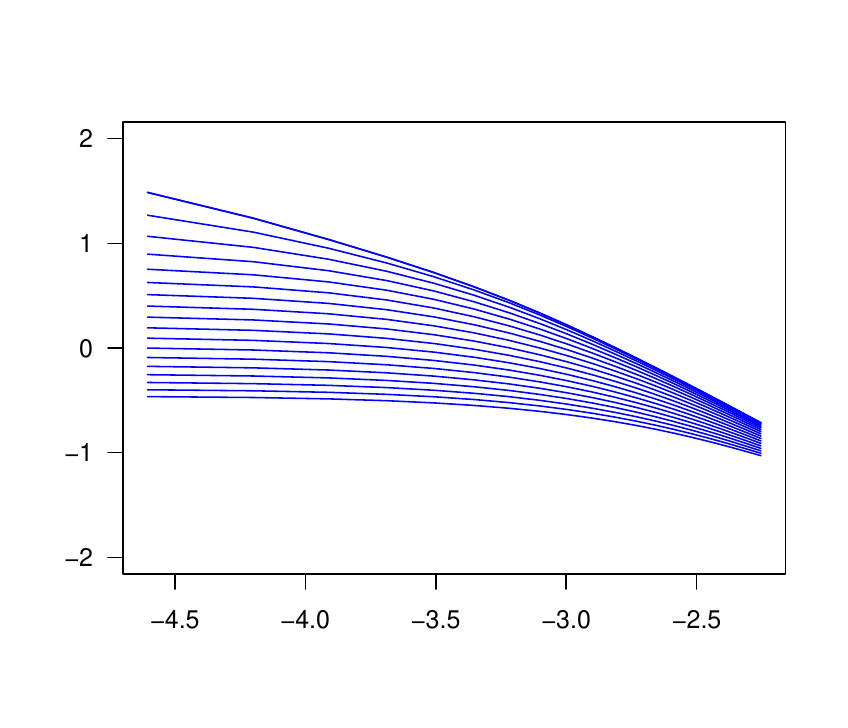}   
    \caption{Plot of $\log(h_2) \mapsto \log\w{\Var}(\w{\pi}_{h_1,h_2}(x^*, y))$. Left graph: for $h_1\in[0.01,0.3]$ with $y= \xi$, right $h_1\in[0.01,0.1]$ with $y=\xi+0.3$.}
    \label{fig:varyxiETvarydiffxi}
\end{figure}

We present here a simple illustration of the results concerning the variance terms when the parameters governing the Hawkes process intensity are known, as analyzed in Propositions \ref{prop: y=} and \ref{prop: y>}. We use the Gaussian kernel to compute the estimator. As shown, the variances depend solely on \(h_2\) (up to a logarithmic term). Specifically, at the point \((x^*, \xi)\), the variance is of order \(1/h_2\). In contrast, at a point \((x^*, y^*)\) where \(y^* > \xi\) (i.e., away from the baseline), the variance is smaller, of order \(1/\sqrt{h_2}\). In Figure \ref{fig:varyxiETvarydiffxi}, we represent the logarithm of the empirical variance at \((x^*, \xi)\) and \((x^*, y^*)\), respectively. For both figures, \(x^*\) is chosen and fixed as the mean value of \(X\). When \(y^*\) is away from \(\xi\), it is set to \(\xi + 0.3\). To keep it simple we investigate here an exponential kernel as described in Equation (\ref{eq: dynamics lambda 1}) and  the set of parameters
$\sigma(\cdot)=0.1$, $a(\cdot)=1$, $b(x)=-6x$,
$(\xi, \alpha, \beta)=(0.5, 0.4, 2)$,  $T=100$, $\Delta = 1/1000$, and we do $1000$ Monte Carlo repetitions. We choose $h_i \in [0.01, 0.3]$ for $i=1, 2$.

Figure \ref{fig:varyxiETvarydiffxi}, left, shows in log-scale the empirical variance of $\w{\pi}_{h_1,h_2}$ for the fixed value of $x^*$ and $y^*= \xi$. It is clear that the slope of the the curve when $y=\xi$ is around $1.5$ for each value of $h_1$ which corresponds to Proposition \ref{prop: y=}. 

When $y^*>\xi$ the behavior is different, as claimed in Proposition \ref{prop: y>}.
For $h_2 \in [0.01, 0.1]$ the slope is around 0.5 as shows right graph of Figure \ref{fig:varyxiETvarydiffxi}.

\begin{figure}[hbtp]
    \centering
    \includegraphics[width=0.4\linewidth]{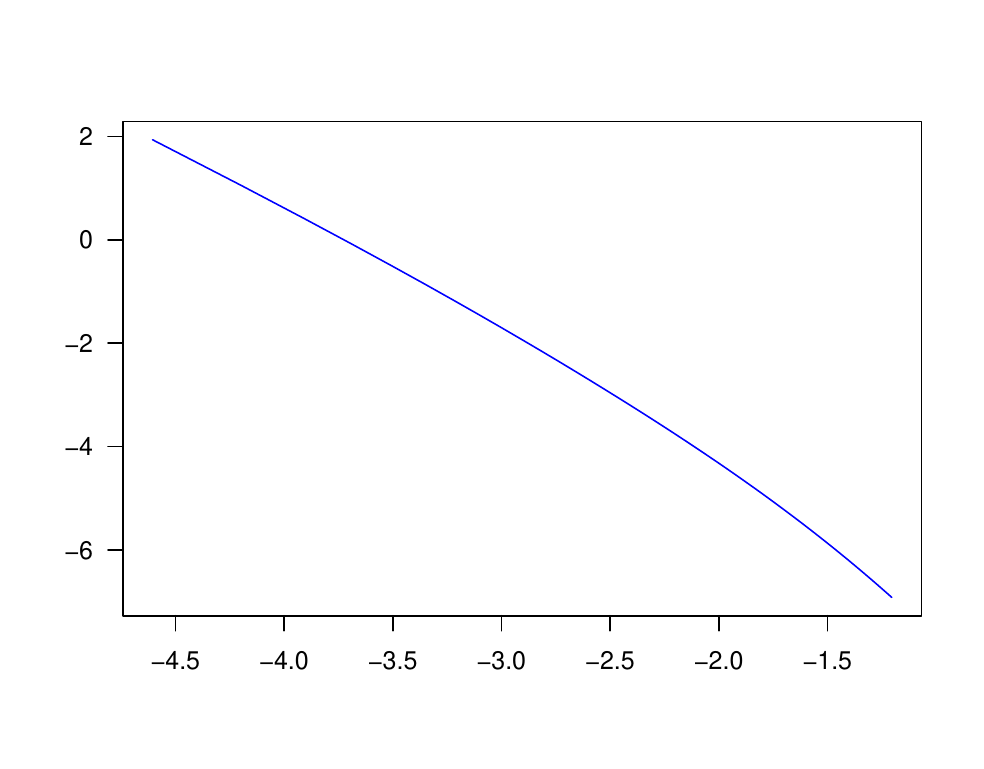}
\includegraphics[width=0.4\linewidth]{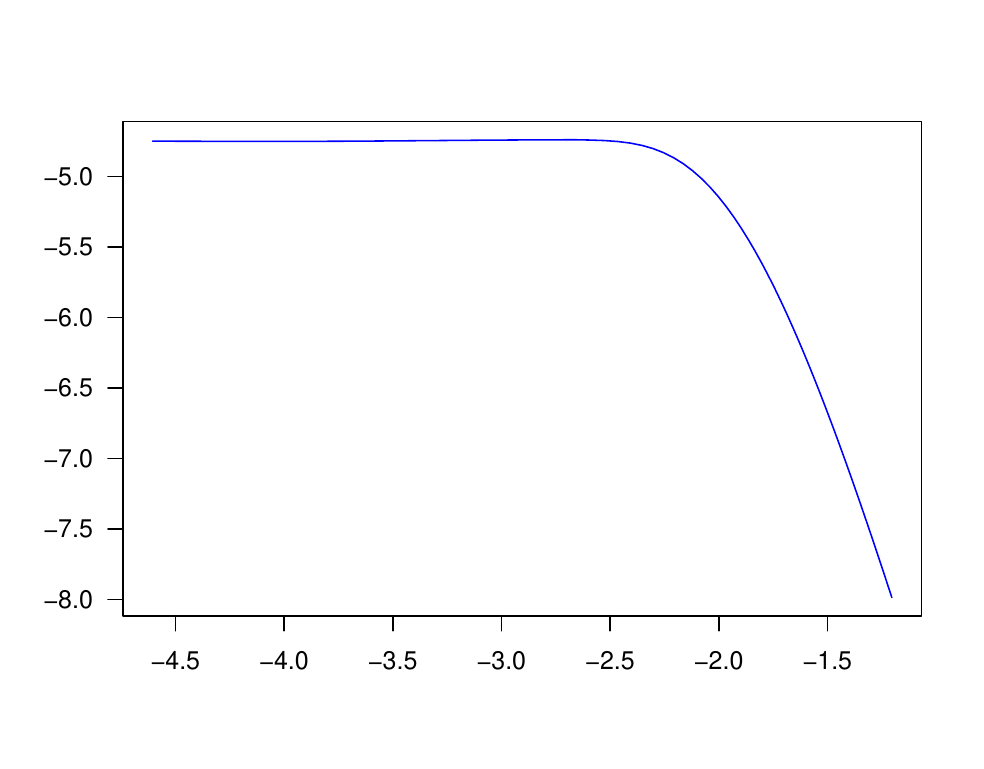}    
    \caption{Plots of $\log(h_2) \mapsto \log\w{\Var}(\w{\pi}_{h_2}(y))$. Left graph: for $y= \xi$, right graph for $y=\xi+0.3$.}
    \label{fig:varyxiandydiffxi}
\end{figure}

We also illustrate the evolution of the variance term in the estimator of the invariant density of the intensity given \(h_2\) (see Proposition \ref{prop: var lambda}). The same simulation setup as described earlier is employed. In this case as well, the variance term depends on whether \(y^* > \xi\) or not. Using a logarithmic scale, the variance remains approximately constant when \(y^* > \xi\), whereas it exhibits a decreasing linear trend in \(h_2\) when \(y^* = \xi\). Figure \ref{fig:varyxiandydiffxi} (left) shows the empirical variance of \(\widehat{\pi}_{h_2}\) in log-scale for a fixed value of \(x^*\) and \(y^* = \xi\), confirming the expected linear decay. For \(y^* > \xi\), the variance initially remains constant up to \(\log(h_2) = -2\), after which it drops sharply. This rapid decline can likely be attributed to the increasing bias term, which causes the variance term to compensate.




\section{Proof of main results}
\label{sec:proofs}
This section is devoted to proving our main results, including the variance bounds presented in Section \ref{sec:estim} and the convergence rates detailed in Corollary \ref{coro:rates}. We then proceed with the proofs of the results discussed in Section \ref{s: estim int lambda} and conclude with the proofs for those in Section \ref{sec:estLambda}.

Let us start by proving the bound on the variance in the case where $y^* = \xi$.

\subsection{Proof of Proposition \ref{prop: y=}}
\begin{proof}
In this proposition, we study the variance of our estimator for \( y^* = \xi \). More precisely, the object of interest is:

\[
\Var(\widehat{\pi}_{h_1, h_2}(x^*,\xi)) = \frac{1}{T^2} \Var\left(\int_0^T \mathbb{K}_{h_1}(x^* - X_t) \mathbb{K}_{h_2}(\xi - \lambda_t) \, \dd t \right),
\]
which can be rewritten as
\[
\frac{1}{T^2} \int_0^T \int_0^T \Cov\left(\mathbb{K}_{h_1}(x^* - X_t) \mathbb{K}_{h_2}(\xi - \lambda_t), \mathbb{K}_{h_1}(x^* - X_s) \mathbb{K}_{h_2}(\xi - \lambda_s)\right) \, \dd t \, \dd s.
\]

Now, consider a function \( f \in L^1(\R^2) \). {We obtain:}

\[
\int_0^T \int_0^T \Cov(f(X_s, \lambda_s), f(X_t, \lambda_t)) \, \dd s \, \dd t = 2 \int_0^T \int_0^s \Cov(f(X_s, \lambda_s), f(X_t, \lambda_t)) \, \dd t \, \dd s,
\]
{and by introducing \( f_c = f - \pi(f) \), the previous inequality  simplifies to}
\[
2 \int_0^T \int_0^s \E[f_c(X_{s-t}, \lambda_{s-t}) f_c(X_0, \lambda_0)] \, \dd t \, \dd s,
\]
a consequence of the stationarity of the pair \( (X_t, \lambda_t) \). By a simple change of variable, we have:
\[
{2}\int_0^T \int_0^s \E[f_c(X_u, \lambda_u) f_c(X_0, \lambda_0)] \, \dd u \, \dd s = 2 \int_0^T (T - u) \E[f_c(X_u, \lambda_u) f_c(X_0, \lambda_0)] \, \dd u.
\]

Let us introduce the following functions for \( f(\cdot, \cdot) = \mathbb{K}_{h_1}(x^* - \cdot) \mathbb{K}_{h_2}(\xi - \cdot) \):
\begin{equation}{\label{eq: def kc}}
k_c(s) := \E[f_c(X_{s}, \lambda_s) f_c(X_0, \lambda_0)]
\end{equation}
and
\begin{equation}{\label{eq: def k}}
k(s) := \E[f(X_{s}, \lambda_s) f(X_0, \lambda_0)].
\end{equation}
Using the above, we can now express the variance as:
\[
\Var(\widehat{\pi}_{h_1, h_2}(x^*,\xi)) \leq \frac{{2}}{T} \int_0^T |k_c(s)| \, \dd s.
\]

To further analyze this integral, we split the time interval \( [0, T] \) into two parts:
\[
[0, T] = [0, D) \cup [D, T],
\]
where \( D{\geq 1} \) will be selected later to obtain the sharpest possible upper bound on the variance. It is straightforward to verify that, because of the boundedness of $\pi$ on the compact support of $f$, \( |k_c(s)| \leq |k(s)|+ C) \). This leads to the inequality:
\begin{equation}{\label{eq: var 19.5}}
\Var(\widehat{\pi}_{h_1, h_2}(x^*, \xi)) \leq \frac{C}{T} \left( \int_0^D {(|k(s)| + 1)} \, \dd s + \int_D^T |k_c(s)| \, \dd s \right).    
\end{equation}

\paragraph{Analysis on \([0, D]\).}
We begin by providing a bound for the first integral above. From the definition of \(k(s)\), conditioning with respect to \((\lambda_u)_{u \in [0,s]}\), we obtain
\begin{equation}{\label{eq: k start I2}}
|k(s)| \leq \E \left[\left| \mathbb{K}_{h_2}(\xi-\lambda_0)\mathbb{K}_{h_2}(\xi-\lambda_s)\right| \E\left[\left| \mathbb{K}_{h_1}(x^*-X_0)\mathbb{K}_{h_1}(x^*-X_s)\right| \bigg| (\lambda_u)_{u \in [0,s]}\right]\right].
\end{equation}

Let \(L_s\) denote the last jump before \(s\) (or {0} if no jump occurs) and \(p_s(x, y)\) the transition density of the process {between $L_s$ and $s$ (i.e. $p_s$ denotes the transition density of the diffusion \eqref{Eq : model diff X intro}, when we have removed all the jumps)}. We will use a bound on the transition density. For unbounded drift, this result is derived from the first point of Theorem 1.2 in \cite{MenPes}. Compared to the bounded drift case, the main difference is the replacement of the starting point by the flow of the starting point in the unbounded drift case. Specifically, for fixed \((s, x) \in \R^+ \times \R\), let \(\theta_{t,s}(x)\) denote the deterministic flow solving \(\dot{\theta}_{t,s}(x)= b(\theta_{t,s}(x))\), \(t \geq 0\), with \(\theta_{s,s}(x) = x\). We have (see also Lemma 1 in \cite{amorino2024minimax}):
\[
p_{t-s}(x, y) \leq \frac{c_1}{\sqrt{t-s}} \exp\left(-c_2 \frac{(\theta_{t,s}(x) - y)^2}{t-s}\right).
\]
Thus, we can write:
\begin{align}{\label{eq: cond X 3.5}}
&\E\left[{|}\mathbb{K}_{h_1}(x^*-X_0)\mathbb{K}_{h_1}(x^*-X_s){|} \bigg| (\lambda_u)_{u \in [0,s]}\right] \nonumber \\
&\leq \E\left[{|}\mathbb{K}_{h_1}(x^*-X_0) \int_{\R} p_{s-L_s}({X_{L_s}}, z) \mathbb{K}_{h_1}(x^*-z) \, \dd z{|} \bigg| (\lambda_u)_{u \in [0,s]}\right] \nonumber \\
&\leq \E\left[{|}\mathbb{K}_{h_1}(x^*-X_0) \int_{\R} \frac{{c_1}}{\sqrt{s-L_s}} \exp\left(-c(z-\theta_{s,L_s}({X_{L_s}}))^2 / (s-L_s)\right) \mathbb{K}_{h_1}(z-x^*) \, \dd z{|} \bigg| (\lambda_u)_{u \in [0,s]}\right] \nonumber \\
&\leq \E\left[{|}\mathbb{K}_{h_1}(x^*-X_0){|}\frac{{c_1}}{\sqrt{s-L_s}} \bigg| (\lambda_u)_{u \in [0,s]}\right] \nonumber \\
&\leq \frac{{c_1}}{\sqrt{s-L_s}} \E\left[{|}\mathbb{K}_{h_1}(x^*-X_0){|} \bigg| (\lambda_u)_{u \in [0,s]}\right],
\end{align}
where \(L_s\) is measurable with respect to \((\lambda_u)_{u \in [0,s]}\). Since \(X_0\) is independent of the \(\sigma\)-algebra generated by \((\lambda_u)_{u \in [0,s]}\), we obtain:
\[
\frac{1}{\sqrt{s-L_s}} \E\left[{|}\mathbb{K}_{h_1}(x^*-X_0){|} \bigg| (\lambda_u)_{u \in [0,s]}\right] \leq \frac{C}{\sqrt{s-L_s}},
\]
using the fact that the \(L^1\)-norm of \(\mathbb{K}_{h_1}\) is simply bounded by a constant {and $\pi$ 
bounded on the support of $\mathbb{K}_{h_1}$}.
We have thus proven that:
\[
|k(s)| \leq C \E\left[\left| \mathbb{K}_{h_2}(\xi-\lambda_0)\mathbb{K}_{h_2}(\xi-\lambda_s)\right| \frac{1}{\sqrt{s-L_s}}\right].
\]

We will now proceed differently depending on whether or not at least one jump occurred in \([0, s]\). Specifically, we introduce the set
\[
\tilde{\Omega}_s := \left\{ \text{There is at least one jump on } [0, s] \right\}.
\]
Thus, we can express
\begin{equation}{\label{eq: k 7}}
|k(s)| \leq C \E\left[\left| \mathbb{K}_{h_2}(\xi-\lambda_0)\mathbb{K}_{h_2}(\xi-\lambda_s)\right| \frac{1}{\sqrt{s - L_s}} \one_{\tilde{\Omega}_s}\right] + C \E\left[\left| \mathbb{K}_{h_2}(\xi-\lambda_0)\mathbb{K}_{h_2}(\xi-\lambda_s)\right| \frac{1}{\sqrt{s - L_s}} \one_{\tilde{\Omega}_s^c}\right].
\end{equation}

Let us first consider the case where at least one jump occurred. In this case, the dynamics of \((\lambda_t)\) after the last jump is deterministic, and we can write {for $s > L_s$,}
\[
\lambda_s = \xi + (\lambda_{L_s} +{\alpha} - \xi)e^{-\beta(s - L_s)}.
\]
Note that the kernel \(\mathbb{K}_{h_2}(\xi - \lambda_s)\) is supported in \([-h_2, h_2]\), meaning that \(\mathbb{K}_{h_2}(\xi - \lambda_s) \neq 0\) only if \(|\lambda_s - \xi| \leq h_2\), {which implies} \(|\lambda_{L_s} +{\alpha} - \xi|e^{-\beta(s - L_s)} \leq h_2\). Since \(|\lambda_{L_s} + {\alpha} - \xi|\) exceeds the jump size {\(\alpha\)}, we have \(e^{-\beta(s - L_s)} \leq \frac{h_2}{{\alpha}}\). As \(h_2 = h_2(T)\) tends to 0 as \(T \to \infty\), we can assume \(\frac{h_2}{\alpha} \leq \frac{1}{2}\). Hence, \((s - L_s) \geq \frac{\log 2}{\beta}\). The first term in \eqref{eq: k 7} is therefore bounded by
\begin{equation}{\label{eq: 22.5}}
C \E\left[\left| \mathbb{K}_{h_2}(\xi - \lambda_0)\mathbb{K}_{h_2}(\xi - \lambda_s)\right|\right] \leq \frac{C}{h_2} \|K\|_\infty \E\left[\left| \mathbb{K}_{h_2}(\xi - \lambda_0)\right|\right] \leq C \frac{1}{h_2},
\end{equation}
where we used the fact that the \(L^1\)-norm of the kernel is bounded by a constant.

Now, let us consider the case where no jumps occurred in \([0, s]\). In this scenario, we have \(L_s = 0\), and the second term in \eqref{eq: k 7} becomes
\[
\frac{C}{\sqrt{s}} \E\left[\left| \mathbb{K}_{h_2}(\xi - \lambda_0)\mathbb{K}_{h_2}(\xi - \lambda_s)\right| \one_{\tilde{\Omega}_s^c}\right].
\]
Since \(\lambda_0\) has a density and the dynamics of \(\lambda_s\) is given by
\[
\lambda_s = \xi + (\lambda_0 - \xi)e^{-\beta s},
\]
substituting this into the expectation above yields
\begin{align}{\label{eq: y=}}
\frac{C}{\sqrt{s}} \E\left[\left| \mathbb{K}_{h_2}(\xi - \lambda_0)\mathbb{K}_{h_2}((\xi - \lambda_0)e^{-\beta s})\right|\right] 
= \frac{C}{\sqrt{s}} \int_{\R} \frac{1}{h_2} K\left(\frac{\xi - z}{h_2}\right) \frac{1}{h_2} K\left(\frac{(z - \xi)e^{-\beta s}}{h_2}\right) \pi_{\lambda_0}(z) \, \dd z,
\end{align}
where \(\pi_{\lambda_0}\) is the density of \(\lambda_0\).
By applying the change of variable \(\frac{\xi - z}{h_2} =: \tilde{z}\), we find that the above expression simplifies to
\begin{equation}{\label{eq: 23.5}}
\frac{C}{\sqrt{s}} \frac{1}{h_2} \int_{\R} K(\tilde{z}) K(\tilde{z}e^{-\beta s}) {\pi_{\lambda_0}}(\xi - \tilde{z}h_2) \, \dd {\tilde{z}} 
\leq \frac{C}{\sqrt{s}} \frac{1}{h_2}.
\end{equation}
Therefore, we conclude that
\begin{equation}{\label{eq: k until D 8}}
\frac{C}{T} \int_{0}^D {(|k(s)| + 1)} \, \dd s \leq \frac{C}{T} \frac{1}{h_2} \int_{0}^D \left(\frac{1}{\sqrt{s}} + 1\right) \, \dd s {+ \frac{C}{T} D} \le \frac{C}{T} \frac{D}{h_2},
\end{equation}
{having remarked that the second term is negligible compared to the first as $h_2 < 1$.}

\paragraph{Analysis of the integral over \([D, T]\).} In this section, we aim to leverage the mixing properties of the process. Specifically, under our hypotheses, the process exhibits exponential beta mixing, as detailed in Theorem 3.8 of \cite{DLL}. {This implies in particular that, for any function \( f_c \) centered under \( \pi \), we have
\begin{equation}{\label{eq: beta mixing}}
 \left\| P_t f_c \right\|_{L^1(\pi)} \le C e^{-\rho t} \left\| f_c \right\|_\infty,   
\end{equation}
where \( (P_t)_{t \ge 0} \) denotes the transition semigroup of the process \( (X, \lambda) \), and \( \rho \) characterizes the exponential \( \beta \)-mixing rate. This definition is consistent with those in \cite{meyn1993stability} and in Definition 1 of \cite{amorino2021invariant}. Bounding $f_c$ with its infinity norm and using beta-mixing inequality as in \eqref{eq: beta mixing}} we obtain:
\[
|k_c(s)| \le \left\| P_s f_c \right\|_{L^1(\pi)} \left\| f_c \right\|_\infty \le C e^{-\rho s} \left\| f_c \right\|_\infty^2.
\] 
Consequently, {replacing $f_c$,} we can establish the following bound on the covariance:
\begin{align*}
|k_c(s)| & \leq C \|\mathbb{K}_{h_1}\|_\infty^2 \|\mathbb{K}_{h_2}\|_\infty^2 e^{-\rho s} \\
& \leq C (h_1 h_2)^{-2} e^{-\rho s}.
\end{align*}
This leads to the estimate
\begin{align}{\label{eq: bound I3}}
\frac{C}{T} \int_{D}^T |k_c(s)| \, \dd s & \leq \frac{C}{T} (h_1 h_2)^{-2} \int_{D}^T e^{-\rho s} \, \dd s \\
& \leq \frac{C}{T} (h_1 h_2)^{-2} e^{-\rho D}. \nonumber
\end{align}
Combining \eqref{eq: k until D 8} and \eqref{eq: bound I3}, we obtain
\begin{align*}
\text{Var}(\w{\pi}_{h_1,h_2}(x^*,\xi)) \leq \frac{C}{T}\left( \frac{D}{h_2} + (h_1 h_2)^{-2} e^{-\rho D }\right).
\end{align*}
To optimize this estimation, we tune the parameter \(D\) such that it becomes as sharp as possible. Notably, by selecting \(D := \max\left(-\frac{2}{\rho}|\log h_1 h_2|, 1\right) \land T\), we can conclude that 
\begin{align*}
\text{Var}(\w{\pi}_{h_1,h_2}(x^*,\xi)) \leq \frac{C}{T}\left(\frac{|\log h_1 h_2|}{h_2} + 1\right) \leq \frac{C}{T}\frac{|\log h_1 h_2|}{h_2},
\end{align*}
as desired.
\end{proof}

We now turn our attention to proving the bound on the variance in the scenario where \(y^* > \xi\). As previously highlighted, this particular case allows for a more favorable convergence rate. The proof of Proposition \ref{prop: y>} sheds light on the underlying reasons for this improved behavior. It is moreover quite involving and relies on the change of measure in the Girsanov style, that we have stated in Section \ref{s:prob}.

\subsection{Proof of Proposition \ref{prop: y>}}
\begin{proof}
In this proposition, we analyze the case where \( y^* > \xi \). The proof of Proposition \ref{prop: y>} follows similar steps to the proof of Proposition \ref{prop: y=}, with the key difference being how we handle \( k(s) \) for \( s \in [0, D] \).

We begin by leveraging the stationarity of the process and splitting the integral into two parts: for \( s \in [0, D) \) and \( [D, T] \). Specifically, we have:
\[
\text{Var}(\w{\pi}_{h_1,h_2}(x^*,y^*)) \leq \frac{C}{T} \left( \int_0^D {(|k(s)| + 1)} \, ds + \int_D^T |k_c(s)| \, ds \right),
\]
with \( k_c(s) \) and \( k(s) \) as defined in \eqref{eq: def kc} and \eqref{eq: def k}, respectively.

Let us first address the second integral, as its analysis mirrors that of the proof of Proposition \ref{prop: y=}. In particular, we can once again invoke the exponential ergodicity of the process, and the bound \eqref{eq: bound I3} still applies. Now, we turn our attention to the first integral, over \( [0, D] \). {In particular, we now aim at proving a bound for $|k(s)|$.}

Observe that \eqref{eq: cond X 3.5} from Proposition \ref{prop: y=} still holds true, as it is independent of \( {y^*} \). We are thus left to bound the expression:
\[
\left| \E\left[\mathbb{K}_{h_2}({y^*} - \lambda_0)\mathbb{K}_{h_2}({y^*} - \lambda_s) \frac{1}{\sqrt{s - L_s}} \right] \right|.
\]

We can rewrite it as:
\[
\E\left[\mathbb{K}_{h_2}(y^* - \lambda_0)\mathbb{K}_{h_2}(y^* - \lambda_s) \frac{1}{\sqrt{s - L_s}} \right]
= \E\left[ \E\left[\mathbb{K}_{h_2}(y^* - \lambda_s) \frac{1}{\sqrt{s - L_s}} \mid \lambda_0 \right] \mathbb{K}_{h_2}(y^* - \lambda_0) \right].
\]

If \( \lambda_0 \notin [y^* - h_2, y^* + h_2] \), then \( \mathbb{K}_{h_2}(y^* - \lambda_0) = 0 \). We thus deduce:
\begin{multline} \label{eq: Kh2 conditione}
\left| \E\left[\mathbb{K}_{h_2}(y^* - \lambda_0)\mathbb{K}_{h_2}(y^* - \lambda_s) \frac{1}{\sqrt{s - L_s}} \right] \right| \\
\leq \int_{y_0 \in [y^* - h_2, y^* + h_2]} \left| \E\left[\mathbb{K}_{h_2}(y^* - \lambda_s) \frac{1}{\sqrt{s - L_s}} \mid \lambda_0 = y_0 \right] \right| \mathbb{K}_{h_2}(y^* - y_0) \, \pi_{\lambda_0}(\dd y_0),
\end{multline}
where \( \pi_{\lambda_0} \) is the stationary distribution of \( \lambda_0 \).

To bound the expression above, we apply the change of measure provided by Theorem \ref{th: girsanov}. For this, some additional notation is necessary. Let us define
\begin{equation} \label{eq: def g q s}
	g^{(q)}(s) := \E\left[ \left| \mathbb{K}_{h_2}(y^* - \lambda_s) \right| \frac{1}{\sqrt{s - L_s}} \one_{\{N_s = q\}} \mid \lambda_0 = y_0 \right], \quad \text{for } q \in \mathbb{N}.
\end{equation}

We now present two lemmas concerning \( g^{(q)}(s) \) that will clarify how to apply Theorem \ref{th: girsanov}. Their proofs are given in Section \ref{s: proof technical}.

\begin{lemme}\label{L: g pour q zero}
	Assume that \( y^* - \xi > 5h_2 \) and \( |y_0 - y^*| \le h_2 \). Then, for \( s \geq \frac{5}{\beta} \frac{h_2}{y^* - \xi} \), we have \( g^{(0)}(s) = 0 \).
\end{lemme} 

\begin{lemme} \label{L: g smaller g hat}
	Let us define \( \ell_g := y^* - \xi - h_2 \), \( \ell_d := y^* - \xi + h_2 \), and the interval \( I = [e^{\beta s} \ell_g + \xi - \lambda_0, e^{\beta s} \ell_d + \xi - \lambda_0] \). For \( q \in \mathbb{N} \setminus \{0\} \), define:
	\begin{equation} \label{eq: def hat g q}
		\hat{g}^{(q)}(s) := \E\left[ \frac{1}{\sqrt{s - T_q}} \one_{\{ c e^{\beta T_1} \left(1 + \sum_{j=2}^q e^{\beta (T_j - T_1)}\right) \in I \}} \one_{\{N_s = q\}} \mid \lambda_0 = y_0 \right].
	\end{equation}
	Then, we have \( g^{(q)}(s) \leq \frac{\|K\|_\infty}{h_2} \hat{g}^{(q)}(s) \) for all \( s > 0 \) and \( q \in \mathbb{N} \setminus \{0\} \).
\end{lemme}

Now, we can state:
\[
\left| \E\left[ \mathbb{K}_{h_2}(y^* - \lambda_s) \frac{1}{\sqrt{s - L_s}} \mid \lambda_0 = y_0 \right] \right| \leq \sum_{q = 0}^\infty g^{(q)}(s).
\]
Since \( y^* > \xi \) is fixed and \( h_2 \) is arbitrarily small, we apply Lemma \ref{L: g pour q zero} for \( s \geq c_1 h_2 \) (with \( c_1 \) given by Lemma \ref{L: g pour q zero}) and deduce that \( g^{(0)}(s) = 0 \). Therefore, we obtain:
\begin{align} \label{eq: before girsanov 42.5}
	\left| \E\left[ \mathbb{K}_{h_2}(y^* - \lambda_s) \frac{1}{\sqrt{s - L_s}} \mid \lambda_0 = y_0 \right] \right| & \leq \sum_{q = 1}^\infty g^{(q)}(s) \nonumber \\
	& \leq \frac{\|K\|_\infty}{h_2} \sum_{q = 1}^\infty \hat{g}^{(q)}(s),
\end{align}
by Lemma \ref{L: g smaller g hat}.

Next, we aim to apply Theorem \ref{th: girsanov} using \( f_1(t) = \frac{1}{\sqrt{t}} \) and \( f_2^{(q)}(t_q, \dots, t_1) = \one_{\{ e^{\beta t_1} \left(1 + \sum_{j=2}^q e^{\beta t_j}\right) \in I \}} \), where \( I \) is defined as in Lemma \ref{L: g smaller g hat}. We need to ensure the assumptions of Theorem \ref{th: girsanov} hold. Clearly, \( f_2 \) is bounded, and we now examine its support.

Recall that the interval \( I \) is given by \( I = [\tilde{\ell}_g, \tilde{\ell}_d] \), with \( \tilde{\ell}_g = e^{\beta s} \ell_g + \xi - \lambda_0 \) and \( \tilde{\ell}_d = e^{\beta s} \ell_d + \xi - \lambda_0 \). Since \( \lambda_0 = y_0 \in (y^* - h_2, y^* + h_2) \), we have:
\begin{align*}
	\tilde{\ell}_g = e^{\beta s} \ell_g + \xi - \lambda_0 &= e^{\beta s}(y^* - \xi - h_2) + \xi - \lambda_0 \\
	& \geq e^{\beta s}(y^* - \xi - h_2) + \xi - y^* - h_2 \\
	& = (e^{\beta s} - 1)(y^* - \xi) - h_2(e^{\beta s} + 1).
\end{align*}
For \( (e^{\beta s} - 1)(y^* - \xi) \geq 4h_2 \), we deduce:
\begin{align} \label{eq: mino lg}
	\tilde{\ell}_g = e^{\beta s} \ell_g + \xi - \lambda_0 &\geq \frac{(e^{\beta s} - 1)(y^* - \xi)}{2} + 2h_2 - h_2e^{\beta s} - h_2 \nonumber \\
	& = (e^{\beta s} - 1)\left( \frac{y^* - \xi}{2} - h_2 \right) \\
	& \geq (e^{\beta s} - 1) \frac{y^* - \xi}{4},
\end{align}
where the last line holds because \( h_2 \leq \frac{y^* - \xi}{4} \). In particular, \( \tilde{\ell}_g > 0 \) under our assumptions.

Next, the constraint \( ce^{\beta u} \left(1 + \sum_{j=1}^{q-1} e^{-\beta t_j}\right) \in [\tilde{\ell}_g, \tilde{\ell}_d] \) leads to \( u \in J^{(q)} \), where \( J^{(q)} \) is the interval:
\[
J^{(q)} = \left[ \frac{1}{\beta} \ln(\tilde{\ell}_g) - \frac{1}{\beta} \ln\left(1 + \sum_{j=1}^{q-1} e^{-\beta t_j}\right) - \frac{\ln(c)}{\beta}, \frac{1}{\beta} \ln(\tilde{\ell}_d) - \frac{1}{\beta} \ln\left(1 + \sum_{j=1}^{q-1} e^{-\beta t_j}\right) - \frac{\ln(c)}{\beta} \right].
\]

In order to apply Theorem \ref{th: girsanov}, we need to find an upper bound, independent on $q$, for the measure $|J^{(q)}|$.
Using the definition of \( J^{(q)} \) and the inequality \( \ln(1 + x) \leq x \), we derive:
\[
|J^{(q)}| = \frac{1}{\beta} \left[\ln(\tilde{l}_d) - \ln(\tilde{l}_g)\right] = \frac{1}{\beta} \ln\left(1 + \frac{\tilde{l}_d - \tilde{l}_g}{\tilde{l}_g}\right)
\]
\[
\leq \frac{\tilde{l}_d - \tilde{l}_g}{\beta \tilde{l}_g} = \frac{e^{\beta s} (l_d - l_g)}{\beta \tilde{l}_g}
\]
\[
\leq \frac{e^{\beta s} (l_d - l_g)}{e^{\beta s} - 1} \cdot \frac{4}{\beta (y^* - \xi)} \quad \text{(by \eqref{eq: mino lg})}
\]
\begin{equation}{\label{eq: 31.5}}
	= \frac{e^{\beta s} 2h_2}{e^{\beta s} - 1} \cdot \frac{4}{\beta (y^* - \xi)}, 
\end{equation}
where, in the last line, we used the definitions of \( \ell_g \) and \( \ell_d \) from Lemma \ref{L: g smaller g hat}.
Hence, we apply Theorem \ref{th: girsanov} with $\overline{J}$ equal to \eqref{eq: 31.5}.
We obtain
\[
\sum_{q = 1}^\infty \hat{g}^{(q)}(s) \leq C_\varepsilon e^{\varepsilon s} \int_0^{\overline{J}} \frac{1}{\sqrt{u}} \, du
 \leq C_\varepsilon e^{\varepsilon s} \sqrt{\overline{J}}.
\]
Replacing \( \overline{J} \) by its expression given in \eqref{eq: 31.5}, we deduce
\[
\sum_{q=1}^{\infty} \hat{g}^{(q)}(s) \leq C_\varepsilon e^{\varepsilon s} \left(\frac{e^{\beta s} h_2}{e^{\beta s} - 1}\right)^{1/2}
\leq C_\varepsilon e^{\varepsilon s} \left(\sqrt{\frac{h_2}{s}} \vee \sqrt{h_2}\right)
\]
for some constant \( C_\varepsilon \) that is independent of \( s \), \( h_2 \), and \( q \). Substituting this upper bound into \eqref{eq: before girsanov 42.5}, we get:

\[
\left|\E\left[\mathbb{K}_{h_2}(y^* - \lambda_s) \frac{1}{\sqrt{s - L_s}} \mid \lambda_0 = y_0\right]\right| \leq C_\varepsilon e^{\varepsilon s} \left( \sqrt{\frac{1}{h_2 s}} \vee \frac{1}{\sqrt{h_2}}\right).
\]

Next, inserting this result into \eqref{eq: Kh2 conditione}, we deduce:
\begin{multline} \label{eq: Kh2 conditione2}
\left| \E\left[\mathbb{K}_{h_2}(y^* - \lambda_0)\mathbb{K}_{h_2}(y^* - \lambda_s) \frac{1}{\sqrt{s - L_s}}\right] \right|  \\
\leq {C_\varepsilon} \left( \frac{1}{\sqrt{h_2 s}} \vee \frac{1}{\sqrt{h_2}} \right) e^{\varepsilon s} \int_{y_0 \in [y^* - h_2, y^* + h_2]} |\mathbb{K}_{h_2}(y^* - y_0)| \pi_{\lambda_0}(\dd y_0) \\
\leq {C_\varepsilon} \left( \frac{1}{\sqrt{h_2 s}} \vee \frac{1}{\sqrt{h_2}} \right) e^{\varepsilon s},
\end{multline}
where we used the fact that \( \pi_{\lambda_0} \) has a bounded density and \( \mathbb{K}_{h_2} \) is bounded by a constant in \( L^1 \).

{Recall that \( D > 1 \), so that we can split the integral over \( [0, D] \) into two parts: \( [0, 1] \) and \( [1, D] \).  
This allows us to clearly identify which of the two terms above attains the maximum.}
It follows that
\begin{align}{\label{eq: first int > 9}}
  \frac{C}{T} \int_0^D {(|k(s)| +1)} \, \mathrm{d}s &\leq \frac{{C_\varepsilon}}{T} \left( \int_0^1 \frac{1}{\sqrt{h_2 s}} e^{s \varepsilon} \, \mathrm{d}s + \int_1^D \frac{1}{\sqrt{h_2}} e^{s \varepsilon} \, \mathrm{d}s {+ D} \right) \nonumber \\
  &\leq \frac{{C_\varepsilon}}{T} \frac{1}{\sqrt{h_2}} \left(1 + e^{c_2 D \varepsilon}\right),
\end{align}
{having used that the last term is negligible compared to the others, as $h_2 < 1$.}

Now, combining \eqref{eq: bound I3} with \eqref{eq: first int > 9}, we obtain the bound:
\[
\text{Var}(\w{\pi}_{h_1,h_2}(x^*, y^*)) \leq \frac{{C_\varepsilon}}{T} \left( \frac{1}{\sqrt{h_2}} \left( 1 + e^{c_2 D \varepsilon} \right) + e^{-\rho D} \frac{1}{(h_1 h_2)^2} \right).
\]

This motivates us to choose \( D := \left[\max\left(- \frac{2}{\rho} \log(h_1 h_2), 1\right) \land T \right] \), so that the bound becomes:
\[
\text{Var}(\w{\pi}_{h_1,h_2}(x^*, y^*)) \leq \frac{{C_\varepsilon}}{T} \left( \frac{1}{\sqrt{h_2}} \left( 1 + (h_1 h_2)^{-\tilde{\varepsilon}} \right) + 1 \right),
\]
which simplifies further to:
\[
\text{Var}(\w{\pi}_{h_1,h_2}(x^*, y^*)) \leq \frac{{C_\varepsilon}}{T} \frac{1}{\sqrt{h_2}} (h_1 h_2)^{-\tilde{\varepsilon}},
\]
where \( \tilde{\varepsilon} \) is a small constant, replacing \( \frac{2 \varepsilon c_2}{\rho} \), as \( \varepsilon \) was 
{fixed} arbitrarily small. Thus, the proof is complete.
\end{proof}

We now proceed to prove the bound on the mean squared error stated in Corollary \ref{coro:rates}. This result follows from the bias-variance decomposition, combined with the variance bounds established for the cases \( y^* = \xi \) and \( y^* > \xi \) in Propositions \ref{prop: y=} and \ref{prop: y>}, respectively.

\subsection{Proof of Corollary \ref{coro:rates}}
\begin{proof} 

To prove this result, we must select the bandwidths \( h_1 \) and \( h_2 \) such that the mean squared error bounds are optimized. 
Let us assume that \( h_1 = h_1(T) = \left(\frac{1}{T}\right)^{a_1} \) and \( h_2 = \left(\frac{1}{T}\right)^{a_2} \) for some \( a_1, a_2 > 0 \) to be determined. In the case where \( y^* = \xi \), we have the following bound for the mean squared error:
\[
c\left(\frac{1}{T}\right)^{2 a_1 \beta_1} + c\left(\frac{1}{T}\right)^{2 a_2 \beta_2} + c\left(\frac{1}{T}\right)^{1 - a_2} \log(T).
\]
This suggests choosing \( a_2 \) such that \( 2 a_2 \beta_2 = 1 - a_2 \), which gives \( a_2 = \frac{1}{2 \beta_2 + 1} \). Meanwhile, \( a_1 \) should be chosen large enough to render the first term negligible compared to the second. Specifically, we should select \( a_1 {\geq} \frac{\beta_2}{\beta_1(2 \beta_2 + 1)} \). 
Thus, for \( y^* = \xi \), we obtain the desired result:
\[
\E\left[\left(\w{\pi}_{h_1^{\rm opt}, h_2^{\rm opt}}(x^*,y^*) - \pi(x^*,y^*)\right)^2\right] \le \frac{c \log(T)}{T^{\frac{2\beta_2}{2 \beta_2 + 1}}},
\]
as required.

Now, let us turn to the case where \( y^* > \xi \). Similarly to the previous scenario, we aim to find \( a_1 \) and \( a_2 \) that minimize the expression:
\[
\left(\frac{1}{T}\right)^{2 a_1 \beta_1} + \left(\frac{1}{T}\right)^{2 a_2 \beta_2} + \left(\frac{1}{T}\right)^{1 - \frac{a_2}{2} - \varepsilon(a_1 + a_2)}.
\]
This leads us to choose \( a_2 \) such that \( 2 a_2 \beta_2 = 1 - \frac{a_2}{2} \), which gives \( a_2 = \frac{2}{4 \beta_2 + 1} \). Similarly to the previous case, we want the first term to be negligible, which implies that \( a_1 {\geq} \frac{2 \beta_2}{\beta_1(4 \beta_2 + 1)} \). Let us choose \( a_1 \) as this value, and redefine
\[
\tilde{\varepsilon} = \varepsilon (a_1 + a_2) = \frac{2 \varepsilon (\beta_1 + \beta_2)}{\beta_1 (4 \beta_2 + 1)}.
\]
Note that since \( \varepsilon > 0 \) is arbitrarily small, \( \tilde{\varepsilon} > 0 \) is also arbitrarily small. Putting everything together, for \( y^* > \xi \), we obtain:
\[
\E\left[\left(\w{\pi}_{h_1^{\rm opt}, h_2^{\rm opt}}(x^*, y^*) - \pi(x^*, y^*)\right)^2\right] \le C \left(\frac{1}{T}\right)^{\frac{4 \beta_2}{4 \beta_2 + 1} - \tilde{\varepsilon}},
\]
which concludes the proof.

\end{proof}

Let us move to the proofs of the results concerning the estimation of the invariant density of the jump intensity $\lambda$, gathered in Proposition \ref{prop: var lambda} and Corollary \ref{c: rate lambda}. 

\subsection{Proof of Proposition \ref{prop: var lambda}}
\begin{proof}
Let us begin by considering the case where \( y^* = \xi \). We closely follow the proof presented in Proposition \ref{prop: y=}. Note that Equation \eqref{eq: var 19.5} still holds, with \( k_c \) and \( k \) defined as in \eqref{eq: def kc} and \eqref{eq: def k}, respectively, but now with \( f(\cdot) = \mathbb{K}_{h_2}(\xi - \cdot) \). For \( s \in [0, D) \), we again split the argument over \( \tilde{\Omega}_s \) and \( \tilde{\Omega}_s^c \), yielding
\[
|k(s)| \leq C \E\left[\left| \mathbb{K}_{h_2}(\xi-\lambda_0)\mathbb{K}_{h_2}(\xi-\lambda_s)\right| \one_{\tilde{\Omega}_s}\right] + C \E\left[\left| \mathbb{K}_{h_2}(\xi-\lambda_0)\mathbb{K}_{h_2}(\xi-\lambda_s)\right| \one_{\tilde{\Omega}_s^c}\right].
\]
This is the same as Equation \eqref{eq: k 7}, but without the extra \( \frac{1}{\sqrt{s - L_s}} \), which arose from the bound on the transition density of the diffusion.

On \( \tilde{\Omega}_s \), we apply Equation \eqref{eq: 22.5}, while on \( \tilde{\Omega}_s^c \), we replace the dynamics of \( \lambda_s \), so that Equations \eqref{eq: y=} and \eqref{eq: 23.5} give
\[
\E\left[\left| \mathbb{K}_{h_2}(\xi-\lambda_0)\mathbb{K}_{h_2}(\xi-\lambda_s)\right| \one_{\tilde{\Omega}_s^c}\right] \le \frac{c}{h_2}.
\]
It follows that
\[
\frac{C}{T} \int_{0}^D |k(s)| \, \dd s \leq \frac{C}{T} \int_{0}^D \frac{1}{h_2} \, \dd s = \frac{C}{T} \frac{D}{h_2}.
\]
On the interval \( [D, T) \), we use exponential ergodicity, as in the proof of Proposition \ref{prop: y=}, which directly provides
\[
\frac{C}{T} \int_{D}^T |k_c(s)| \, \dd s \leq \frac{C}{T} h_2^{-2} e^{-\rho D},
\]
leading us to choose \( D := \max\left(-\frac{2}{\rho}\log|h_2|, 1\right) \land T \). This yields
\[
\Var(\w{\pi}_{h_2}(\xi)) \leq \frac{C}{T} \frac{\log|h_2|}{h_2},
\]
as desired.

Let us now move to the case where \( y^* > \xi \). In this case, to avoid integrability problems near \( 0 \), we split the integral into three parts, leading to
\[
\Var(\widehat{\pi}_{h_2}(y^*)) \leq \frac{C}{T} \left( \int_0^\delta |k(s)| \, \dd s + \int_\delta^D |k(s)| \, \dd s + \int_D^T |k_c(s)| \, \dd s \right),
\]
with \( \delta \) and \( D \) to be chosen later.

For \( s \in [0, \delta) \), using the Cauchy-Schwarz inequality and the stationarity of the process, we get
\[
|k(s)| \le \Var(\mathbb{K}_{h_2}(y^* - \lambda_0))^\frac{1}{2} \Var(\mathbb{K}_{h_2}(y^* - \lambda_s))^\frac{1}{2} = \Var(\mathbb{K}_{h_2}(y^* - \lambda_0)) \le \int_{\R} (\mathbb{K}_{h_2}(y^* - z))^2 \pi_{\lambda_0}(z) \, dz,
\]
where \( \pi_{\lambda_0} \) is the density of \( \lambda_0 \). By bounding the infinity norm of \( \mathbb{K}_{h_2} \), we get \( |k(s)| \le \frac{c}{h_2^2} \), so that
\[
\frac{c}{T} \int_0^\delta |k(s)| \, \dd s \le \frac{c \delta}{T h_2^2}.
\]

Next, for \( s \in [\delta, D) \), we follow the approach in Proposition \ref{prop: y>}, noting the absence of the term \( \frac{1}{\sqrt{s - L_s}} \) that originated from the bound on the transition density of \( X \). In analogy with the definitions of \( g^{(q)} \) and \( \hat{g}^{(q)} \), we introduce
\[
g_\lambda^{(q)}(s) := \E\left[ \left| \mathbb{K}_{h_2}(y^* - \lambda_s) \right| \one_{\{N_s = q\}} \mid \lambda_0 = y_0 \right], \quad \text{for } q \in \mathbb{N},
\]
and
\[
\hat{g}_\lambda^{(q)}(s) := \E\left[ \one_{\{ \alpha e^{\beta T_1} \left(1 + \sum_{j=2}^q e^{\beta (T_j - T_1)}\right) \in I \}} \one_{\{N_s = q\}} \mid \lambda_0 = y \right].
\]
It is straightforward to verify that the results of Lemmas \ref{L: g pour q zero} and \ref{L: g smaller g hat} still hold.

We then apply Theorem \ref{th: girsanov} with \( \tilde{f}_1(t) = 1 \) and \( f_2^{(q)}(t) \) as in the proof of Proposition \ref{prop: y>}, confirming that the assumptions on \( f_2^{(q)} \) required for Theorem \ref{th: girsanov} are satisfied, as they were verified in Proposition \ref{prop: y>}. This gives
\[
\sum_{q = 1}^\infty \hat{g}_\lambda^{(q)}(s) \le C_\varepsilon e^{\varepsilon s} \overline{J} \le C_\varepsilon e^{\varepsilon s} \left( \frac{h_2}{s} \lor h_2 \right),
\]
as shown in Equation \eqref{eq: 31.5}. Therefore, similarly as in \eqref{eq: Kh2 conditione2}, we have
\[
\left| \E\left[\mathbb{K}_{h_2}(y^* - \lambda_0)\mathbb{K}_{h_2}(y^* - \lambda_s) \right] \right| \leq C_\varepsilon e^{\varepsilon s} \left( \frac{1}{s} \lor 1 \right).
\]
It follows that
\[
\frac{C}{T} \int_\delta^D |k(s)| \, \dd s \le \frac{C}{T} \int_\delta^1 \frac{1}{s} e^{\varepsilon s} \, \dd s + \frac{C}{T} \int_1^D e^{\varepsilon s} \, \dd s \le \frac{C}{T}(\log \delta + e^{\varepsilon D}).
\]

Finally, we use exponential ergodicity on the interval \( [D, T) \). Combining all the pieces, we obtain
\[
\text{Var}(\w{\pi}_{h_2}(\xi)) \leq \frac{C}{T} \left( \frac{\delta}{h_2^2} + \log \delta + e^{\varepsilon D} + h_2^{-2} e^{- \rho D} \right).
\]
This leads to choosing \( \delta = h_2^2 \) and \( D := \max\left(-\frac{2}{\rho}|\log h_2|, 1\right) \land T \), yielding 
\[
\text{Var}(\w{\pi}_{h_2}(\xi)) \leq \frac{C}{T} \left( |\log h_2| + h_2^{- \varepsilon} \right),
\]
as desired.
\end{proof}

\subsection{Proof of Corollary \ref{c: rate lambda}}
\begin{proof}
The proof is again based on a bias-variance decomposition. Following the same steps as in Corollary \ref{coro:rates}, since the bound on the variance in Proposition \ref{prop: var lambda} is of the same order as that in Proposition \ref{prop: y=}, we obtain the same optimal rate for the bandwidth \( h_2 \). Consequently, this leads to the same convergence rate, \( \frac{c \log(T)}{T^{\frac{2\beta_2}{2 \beta_2 + 1}}} \) when $y^*= \xi$.

Now, let us analyze the case where \( y^* > \xi \). Here, we seek the bandwidth \( h_2(T) = \left(\frac{1}{T}\right)^{a_2} \) that minimizes
\[
h_2^{2 \beta_2} + \frac{C}{T}(|\log h_2| + h_2^{- \varepsilon}).
\]
It is clear that, by choosing \( a_2 \) small enough to satisfy \( 2 a_2 \beta_2 < 1 - a_2 \varepsilon \), the bias term becomes negligible compared to the variance. Thus, for arbitrarily small \( \varepsilon > 0 \), the convergence rate is
\[
\frac{C}{T}(\log T + T^{\varepsilon}),
\]
which completes the proof.

\end{proof}

We now proceed to the proofs of the results presented in Section \ref{sec:estLambda}, focusing on the case where the jump intensity is no longer assumed to be known, and must first be estimated. In the following subsection, we demonstrate that the estimator we propose attains parametric convergence rates.

\subsection{Proof of Theorem \ref{th: bound lambda hat}}
\begin{proof}
Observe that, from the dynamics of the intensity \(\lambda\) outlined in \eqref{eq: dynamics lambda 1}, we can express the dynamics of the estimator \(\hat{\lambda}\) introduced in \eqref{eq: def lambda hat 2} as follows:
\[
\hat{\lambda}_t= \hat{\xi}+ \int_{0}^{{t-}} \hat{\alpha} e^{-\hat{\beta}(t-u)} \dd N_u + (\tilde{\lambda}_0-\hat{\xi}) e^{-\hat{\beta} t}.
\]
This leads to:
\begin{align*}
\lambda_t - \hat{\lambda}_t & = \xi - \hat{\xi} + \int_{0}^{{t-}} \left(\alpha e^{-\beta(t-u)} - \hat{\alpha} e^{-\hat{\beta}(t-u)}\right) \dd N_u + (\lambda_0 - \xi)e^{-\beta t} - (\tilde{\lambda}_0 - \hat{\xi}) e^{-\hat{\beta} t} 
\end{align*}
and we denote 
$$M_t :=\int_{0}^{{t-}} \left(\alpha e^{-\beta(t-u)} - \hat{\alpha} e^{-\hat{\beta}(t-u)}\right) \dd N_u $$
We have required that both \(\beta\) and \(\hat{\beta}\) are lower bounded by some \(\tilde{\varepsilon}\). This implies:
\begin{align*}
|\lambda_t - \hat{\lambda}_t| & \le |\xi - \hat{\xi}|(1 + e^{-\tilde{\varepsilon}t}) + |M_t| + |\lambda_0 - \tilde{\lambda}_0| e^{-\tilde{\varepsilon}t}.
\end{align*}
From here, we obtain:
\begin{align}{\label{eq: norm q 2.5}}
\E\left[|\lambda_t - \hat{\lambda}_t|^q\right] & \le c \E[|\xi - \hat{\xi}|^q] + c \E[|M_t|^q] + c \E[|\lambda_0 - \tilde{\lambda}_0|^q] e^{-q \tilde{\varepsilon} t} \nonumber \\
& \le \frac{c}{T^{\frac{q}{2}}} + c \E[|M_t|^q] + c \E[|\lambda_0 - \tilde{\lambda}_0|^q] e^{-q \tilde{\varepsilon} t}, 
\end{align}
where we have utilized \eqref{eq: parametric rate 1.5} in the last inequality. 

Next, let us focus on the \(q\)-moments of \(M_t\). To evaluate these moments, we will employ Kunita's inequality. We refer to the Appendix of \cite{JacPro} for a proof of this inequality in a general framework, and on page 52 of \cite{amorino2020contrast}, a proof is provided in a context closer to ours. 

Recall that for a compensated Poisson random measure \(\tilde{\mu} = \mu - \bar{\mu}\) and a jump coefficient \(j(x, z)\), Kunita's inequality asserts that, for any \(q \ge 2\):
\[
\E\left[\left|\int_0^t \int_{\mathbb{R}} j(X_{s^{-}}, z) \tilde{\mu}(ds, dz)\right|^q\right] \le C \E\left[\int_0^t \int_{\mathbb{R}} |j(X_{s^{-}}, z)|^q \bar{\mu}(\dd s, \dd z)\right] + C\E\left[\left|\int_0^t \int_{\mathbb{R}} j^2(X_{s^{-}}, z) \bar{\mu}(\dd s, \dd z)\right|^{\frac{q}{2}}\right],
\]
where \(\bar{\mu}\) is the compensator. 

In our case, we will apply Kunita's inequality to the compensated measure \(\dd \tilde{N}_u\) associated with the (non-compensated) measure \(\dd N_u\) appearing in \(M_t\). Here, the compensator is given by \(\lambda_u \, du\). 
Thus, for any \(q \ge 2\), we have:
\begin{align}{\label{eq: bound Mt 3}}
\E\left[|M_t|^q]\right] & \le C \int_0^t \E\left[|\alpha e^{-\beta(t-u)} - \hat{\alpha} e^{-\hat{\beta}(t-u)}|^q \lambda_u\right] \dd u \nonumber\\
& \quad + C  \E\left[\left( \int_0^t |\alpha e^{-\beta(t-u)} - \hat{\alpha} e^{-\hat{\beta}(t-u)}|^2 \lambda_u \, \dd u\right)^{\frac{q}{2}}\right] \nonumber \\
& \quad + C \E\left[\left( \int_0^t |\alpha e^{-\beta(t-u)} - \hat{\alpha} e^{-\hat{\beta}(t-u)}| \lambda_u \, \dd u\right)^q\right] \nonumber \\
& =: M_1 + M_2 + M_3. 
\end{align}

Let us begin by analyzing \(M_1\). By strategically adding and removing \(\alpha e^{-\hat{\beta}(t-u)}\) within the integral, we obtain the following estimate:
\[
M_1 \le C \int_0^t \E[|\alpha - \hat{\alpha}|^q e^{-\hat{{\beta}}(t-u)q} \lambda_u] \, \dd u + C \int_0^t |\alpha|^q \E[|e^{-\beta(t-u)} - e^{-\hat{\beta}(t-u)}|^q \lambda_u] \, \dd u.
\]
Since \(\hat{\alpha}\), \(\hat{\xi}\), and \(\hat{\beta}\) are computed from a different dataset, they are independent of \(\lambda_u\). Moreover, for the first term, we use the fact that \(\hat{\beta} > \tilde{\varepsilon}\), and in the second term, we apply the mean value {theorem}, recalling that \(\alpha\) is bounded. This yields:
\[
M_1 \le C \int_0^t \E[|\alpha - \hat{\alpha}|^q] e^{-\tilde{\varepsilon}(t-u)q} \E[\lambda_u] \, \dd u + C \int_0^t \E[|\hat{\beta} - \beta|^q (t-u)^q e^{-\tilde{\beta}(t-u){q}}] \E[\lambda_u] \, \dd u,
\]
for some \(\tilde{\beta} \in [\beta, \hat{\beta}]\). Under our hypotheses, {the Hawkes process is stationary and linear, so that $\mathbb{E}\left[\lambda_u\right]=\xi /(1-\alpha / \beta)$.} Consequently, using \eqref{eq: parametric rate 1.5}, we find:
\begin{equation}{\label{eq: M1 4}}
M_1 \le \frac{C}{T^{\frac{q}{2}}} \int_0^t e^{-\tilde{\varepsilon}(t-u)q} \, \dd u + \frac{C}{T^{\frac{q}{2}}} \int_0^t (t-u)^{{q}} e^{-\tilde{\varepsilon}(t-u)q} \, \dd u \le \frac{C}{T^{\frac{q}{2}}},
\end{equation}
as desired. 

Next, we analyze \(M_2\). We can express it as:
\[
M_2 = C \E\left[\left( \int_0^t |\alpha e^{-(\beta - \frac{\tilde{\varepsilon}}{2})(t-u)} - \hat{\alpha} e^{-(\hat{\beta} - \frac{\tilde{\varepsilon}}{2})(t-u)}|^2 e^{-\tilde{\varepsilon}(t-u)} \lambda_u \, \dd u\right)^{\frac{q}{2}} \right].
\]
{Consider the measure \( \mu(\dd s) := e^{- \tilde{\varepsilon}(t - s)} \dd s \).  
Possibly, 
\( \int_0^t \mu(\dd s) \neq 1 \), which leads us to apply a 
re-normalized version of Jensen's inequality, namely:
\[
\left( \frac{ \int_0^t f(s) \, \mu(\dd s) }{ \int_0^t \mu(\dd s) } \right)^{\frac{q}{2}} 
\le \frac{1}{ \int_0^t \mu(\dd s) } \int_0^t f(s)^{\frac{q}{2}} \, \mu(\dd s).
\]
This motivates multiplying and dividing by \( \left( \int_0^t e^{- \tilde{\varepsilon}(t - s)} \, \dd s \right)^{\frac{q}{2}} \), yielding:}

\begin{align*}
M_2 & \le C \left( \int_0^t e^{-\tilde{\varepsilon}(t-u)} \, \dd u\right)^{\frac{q}{2} - 1} \E\left[\int_0^t |\alpha e^{-(\beta - \frac{\tilde{\varepsilon}}{2})(t-u)} - \hat{\alpha} e^{-(\hat{\beta} - \frac{\tilde{\varepsilon}}{2})(t-u)}|^q \lambda_u^{\frac{q}{2}} e^{-\tilde{\varepsilon}(t-u)} \, \dd u\right] \\
& \le C \int_0^t \E[|\alpha e^{-(\beta - \frac{\tilde{\varepsilon}}{2})(t-u)} - \hat{\alpha} e^{-(\hat{\beta} - \frac{\tilde{\varepsilon}}{2})(t-u)}|^q \lambda_u^{\frac{q}{2}}] \, \dd u.
\end{align*}
Utilizing the same reasoning as for \(M_1\), we note that under our 
{hypothesis,
by Proposition 1 in \cite{ADBGL}
for any \(q \ge 2\), \(\sup_{u} \E[\lambda_u^{\frac{q}{2}}] < \infty\)} and that {\(\min(\beta - \frac{\tilde{\varepsilon}}{2}, \hat{\beta} - \frac{\tilde{\varepsilon}}{2}) > \frac{\tilde{\varepsilon}}{2}\)}. This leads us to conclude:
\begin{equation}{\label{eq: M2 5}}
M_2 \le \frac{C}{T^{\frac{q}{2}}}.
\end{equation}
Now, we proceed to apply Jensen's inequality to \(M_3\). We write:
\begin{align*}
M_3 & = C \E\left[\left( \int_0^t |c e^{-(\alpha - \frac{\tilde{\varepsilon}}{2})(t-u)} - \hat{c} e^{-(\hat{\alpha} - \frac{\tilde{\varepsilon}}{2})(t-u)}| e^{-\frac{\tilde{\varepsilon}}{2}(t-u)} \lambda_u \, \dd u\right)^q \right] \\
& \le C \left( \int_0^t e^{-\frac{\tilde{\varepsilon}}{2}(t-u)} \, \dd u\right)^{q-1} \E\left[\int_0^t |c e^{-(\alpha - \frac{\tilde{\varepsilon}}{2})(t-u)} - \hat{c} e^{-(\hat{\alpha} - \frac{\tilde{\varepsilon}}{2})(t-u)}|^q \lambda_u^{q} e^{-\frac{\tilde{\varepsilon}}{2}(t-u)} \,\dd u\right].
\end{align*}
Again, the first integral is bounded, and since \(\sup_u \E[\lambda_u^q] < \infty\) for any \(q \ge 2\), we can apply the arguments used for \(M_1\) to \(M_3\) as well. This leads us to conclude:
\begin{equation}{\label{eq: M3 6}}
M_3 \le \frac{C}{T^{\frac{q}{2}}}.
\end{equation}
Substituting \eqref{eq: M1 4}, \eqref{eq: M2 5}, and \eqref{eq: M3 6} back into \eqref{eq: bound Mt 3}, we obtain \(\E[|M_t|^q] \le \frac{C}{T^{\frac{q}{2}}}\). Together with \eqref{eq: norm q 2.5}, this concludes the proof of the theorem.
\end{proof}

The bound on the estimation error of \(\lambda\) established above is instrumental in controlling the error incurred when transitioning from the estimator \(\tilde{\pi}_{h_1, h_2}\), which is based solely on the observed process \(X\), to \(\hat{\pi}_{h_1, h_2}\), which can be utilized when both the pair \((X, \lambda)\) is observed. This is demonstrated in detail in the proof of Proposition \ref{prop: Ta} below.

\subsection{Proof of Proposition \ref{prop: Ta}}
\begin{proof}
From the definitions of \(\tilde{\pi}_{h_1, h_2}\) and \(\hat{\pi}_{h_1, h_2}\), and applying Jensen's inequality, we find that the error we aim to analyze is bounded as follows:

\begin{align*}
& \frac{1}{(T - \Tm)^2} \E\left[\left(\int_{\Tm}^T  \mathbb{K}_{h_1}(x^* - X_u)(\mathbb{K}_{h_2}(y^* - \hat{\lambda}_u)- \mathbb{K}_{h_2}(y^* - {\lambda}_u)) \dd u\right)^2\right] \\
& \le \frac{C}{T - \Tm} \int_{\Tm}^T \E[\mathbb{K}_{h_1}^2(x^* - X_u)(\mathbb{K}_{h_2}(y^* - \hat{\lambda}_u) - \mathbb{K}_{h_2}(y^* - {\lambda}_u))^2] \dd u.
\end{align*}
Using the mean value theorem, we have:
\begin{align*}
|\mathbb{K}_{h_2}(y^* - \hat{\lambda}_u) - \mathbb{K}_{h_2}(y^* - {\lambda}_u)| & = \frac{1}{h_2} \left|K\left(\frac{y^* - \hat{\lambda}_u}{h_2}\right) - K\left(\frac{y^* - {\lambda}_u}{h_2}\right)\right| \\
& \le \frac{1}{h_2} \left\| K' \right\|_\infty \frac{|\hat{\lambda}_u - {\lambda}_u|}{h_2}.
\end{align*}

This leads us to the following bound:

\begin{align*}
& \frac{C}{T - \Tm} \frac{1}{h_2^4} \int_{\Tm}^T \E\left[\mathbb{K}_{h_1}^2(x^* - X_u)|\hat{\lambda}_u - {\lambda}_u|^2\right] \dd u \\
& \le \frac{C}{T - \Tm} \frac{1}{h_2^4} \int_{\Tm}^T \E[|\mathbb{K}_{h_1}(x^* - X_u)|^{2p}]^{\frac{1}{p}} \E\left[|\hat{\lambda}_u - {\lambda}_u|^{2q}\right]^{\frac{1}{q}} \dd u \\
& \le \frac{C}{T - \Tm} \frac{1}{h_2^4} \frac{1}{T} \int_{\Tm}^T \E\left[|\mathbb{K}_{h_1}(x^* - X_u)|^{2p}\right]^{\frac{1}{p}} \dd u,
\end{align*}

where we have applied Hölder's inequality and Theorem \ref{th: bound lambda hat}. Note that since the \(L^1\) norm of \(\mathbb{K}_{h_1}\) is bounded by a constant, and its supremum norm is bounded by \(\frac{C}{h_1}\), we obtain:
\[
\E\left[|\mathbb{K}_{h_1}(x^* - X_u)|^{2p}\right] = \E[|\mathbb{K}_{h_1}(x^* - X_u)|^{2p-1} |\mathbb{K}_{h_1}(x^* - X_u)|] \le \left(\frac{C}{h_1}\right)^{2p - 1}.
\]
This motivates us to choose \(p\) close to \(1\). Thus, for an arbitrarily small \(\varepsilon > 0\), we can set \(p = 1 + \varepsilon/{2}\). This leads to the conclusion:
\begin{equation}{\label{eq: end pi hat}}
\E\left[(\tilde{\pi}_{h_1, h_2}(x^*, y^*) - \hat{\pi}_{h_1, h_2}(x^*, y^*))^2\right] \le \frac{C}{T} \frac{1}{h_2^4} \frac{1}{{h_1}^{1 + \varepsilon}},
\end{equation}
as desired. This concludes the first part of the proof. The second part focuses on the selection of rate-optimal bandwidths, \(h_1\) and \(h_2\). Indeed, using \eqref{eq:decompHatlambda}, along with the results from Section \ref{sec:estim}, we obtain the bound:
\begin{equation}{\label{eq: decomp bound}}
\E\left[(\tilde{\pi}_{h_1, h_2}(x^*, y^*) - \pi(x^*, y^*))^2\right] \le C(h_1^{2\beta_1} + h_2^{2{\beta_2}}) + \frac{C}{T} \frac{1}{h_2^4} \frac{1}{{h_1}^{1 + \varepsilon}},
\end{equation}
where we have already accounted for the fact that, since \(h_1, h_2 \leq 1\), the variance bounds derived from Propositions \ref{prop: y=} and \ref{prop: y>} are smaller than those in \eqref{eq: end pi hat}.

Next, as in the proof of Corollary \ref{coro:rates}, we seek the bandwidths \(h_1(T)= \left(\frac{1}{T}\right)^{a_1}\) and \(h_2(T)= \left(\frac{1}{T}\right)^{a_2}\) that achieve an optimal trade-off. This leads us to the following system of equations:

\[
\begin{cases}
    2a_1 \beta_1 = 2a_2 \beta_2, \\
    2a_1 \beta_1 = 1 - (1 + \varepsilon)a_1 - 4a_2.
\end{cases}
\]

Solving this system gives the two optimal bandwidths:

\begin{equation}{\label{eq: h lambdahat}}
 h_1(T) = \left(\frac{1}{T}\right)^{\frac{\beta_2}{\beta_1} \frac{1}{\beta_2(2 + \frac{1 + \varepsilon}{\beta_1}) + 4}}, \quad h_2(T) = \left(\frac{1}{T}\right)^{\frac{1}{\beta_2(2 + \frac{1 + \varepsilon}{\beta_1}) + 4}}.   
\end{equation}
Substituting these into \eqref{eq: decomp bound} yields:

\[
\E\left[(\tilde{\pi}_{h_1, h_2}(x^*, y^*) - \pi(x^*, y^*))^2\right] \le \left(\frac{1}{T}\right)^{\frac{2 \beta_2}{\beta_2(2 + \frac{1 + \varepsilon}{\beta_1}) + 4}}.
\]

This concludes the proof once we redefine \(\tilde{\varepsilon}\) such that:

\[
\frac{2 \beta_2}{\beta_2(2 + \frac{1}{\beta_1}) + 4 + \varepsilon\frac{\beta_2}{\beta_1}} = \frac{2 \beta_2}{\beta_2(2 + \frac{1}{\beta_1}) + 4} - \tilde{\varepsilon},
\]

and note that \(\tilde{\varepsilon}\) can be made arbitrarily small, just as \(\varepsilon\) was.

\end{proof}


\section{Proof of technical results}{\label{s: proof technical}}

This section is dedicated to proving the technical results stated earlier, which are essential for establishing our main findings but have not yet been demonstrated. We begin by presenting the proof of the bound on the bias term, as outlined in Proposition \ref{prop:biais}.

\subsection{Proof of Proposition \ref{prop:biais}}
\begin{proof}

The proof closely follows the reasoning presented in Proposition 2 of \cite{amorino2021invariant}. Starting from the definition of \(\w{\pi}(x^*,y^*)\), we can express

$$\E[\w{\pi}(x^*,y^*)] = \int_{\R^2} \mathbb{K}_{h_1}(x^* - u)\mathbb{K}_{h_2}(y^* - v) \pi(u,v) \dd u \dd v. $$
Thus, we have

\begin{align}{\label{eq: bias 10}}
 |\E[\w{\pi}(x^*,y^*)] - \pi(x^*,y^*)| &= \left|\frac{1}{h_1 h_2} \int_{\R^2} K\left(\frac{x^* - u}{h_1}\right) K\left(\frac{y^* - v}{h_2}\right) \pi(u,v) \dd u \dd v - \pi(x^*,y^*)\right| \nonumber \\
 &= \left|\int_{\R^2} K(\tilde{u}) K(\tilde{v}) \pi(x^* - h_1 \tilde{u}, y^*-h_2 \tilde{v}) \dd\tilde{u} \dd\tilde{v} - \pi(x^*,y^*)\right| \nonumber \\
&= \left|\int_{\R^2} K(\tilde{u}) K(\tilde{v}) \left[\pi(x^* - h_1 \tilde{u}, y^* - h_2 \tilde{v}) - \pi(x^*,y^*)\right] \dd\tilde{u} \dd\tilde{v}\right|,
\end{align}
where we applied the change of variables \(\tilde{u} := \frac{x^* - u}{h_1}\) and \(\tilde{v} := \frac{y^* - v}{h_2}\), and used the fact that \(\int_{\R^2} K(\tilde{u}) K(\tilde{v}) \dd \tilde{u} \dd \tilde{v} = 1\) by the definition of the kernel function.
We then apply Taylor's formula to the partial functions \( t \mapsto \pi(x^* - h_1 \tilde{u}, t) \) and \( t \mapsto \pi(t, y^*) \), yielding

\[
\pi(x^* - h_1 \tilde{u}, y^* - h_2 \tilde{v}) = \pi(x^*, y^*) + \sum_{k = 1}^{\lfloor \beta_2 \rfloor - 1} \frac{\partial_y^k \pi(x^* - h_1 \tilde{u}, y^*)}{k!} ({-}h_2 \tilde{v})^k + \frac{\partial_y^{\lfloor \beta_2 \rfloor} \pi(x^* - h_1 \tilde{u}, y^* - h_2 \tau_2 \tilde{v})}{\lfloor \beta_2 \rfloor!} ({-}h_2 \tilde{v})^{\lfloor \beta_2 \rfloor} 
\]
\[
+ \sum_{k = 1}^{\lfloor \beta_1 \rfloor - 1} \frac{\partial_x^k \pi(x^*, y^*)}{k!} ({-}h_1 \tilde{u})^k + \frac{\partial_x^{\lfloor \beta_1 \rfloor} \pi(x^* - \tau_1 h_1 \tilde{u}, y^*)}{\lfloor \beta_1 \rfloor!} ({-}h_1 \tilde{u})^{\lfloor \beta_1 \rfloor},
\]
where \(0 \le \tau_1, \tau_2 \le 1\). Substituting this expansion into \eqref{eq: bias 10}, and recalling that the kernel is of order \( M = \max(\lfloor \beta_1 \rfloor, \lfloor \beta_2 \rfloor) \), we obtain

\[
\int_{\R^2} K(\tilde{u}) K(\tilde{v}) \left( \sum_{k = 1}^{\lfloor \beta_2 \rfloor - 1} \frac{\partial_y^k \pi(x^* - h_1 \tilde{u}, y^*)}{k!} ({-}h_2 \tilde{v})^k + \sum_{k = 1}^{\lfloor \beta_1 \rfloor - 1} \frac{\partial_x^k \pi(x^*, y^*)}{k!} ({-}h_1 \tilde{u})^k  \right) \dd \tilde{u} \dd \tilde{v} = 0.
\]

Therefore, \eqref{eq: bias 10} becomes

\begin{align}{\label{eq: order 11}}
& \int_{\R^2} K(\tilde{u}) K(\tilde{v}) \left( \frac{\partial_y^{\lfloor \beta_2 \rfloor} \pi(x^* - h_1 \tilde{u}, y^* - h_2 \tau_2 \tilde{v})}{\lfloor \beta_2 \rfloor!} ({-}h_2 \tilde{v})^{\lfloor \beta_2 \rfloor} + \frac{\partial_x^{\lfloor \beta_1 \rfloor} \pi(x^* - \tau_1 h_1 \tilde{u}, y^*)}{\lfloor \beta_1 \rfloor!} ({-}h_1 \tilde{u})^{\lfloor \beta_1 \rfloor} \right) \dd \tilde{u} \dd \tilde{v} \nonumber \\
= & \int_{\R^2} K(\tilde{u}) K(\tilde{v}) \Bigg( \frac{({-}h_2 \tilde{v})^{\lfloor \beta_2 \rfloor}}{\lfloor \beta_2 \rfloor!} \left(\partial_y^{\lfloor \beta_2 \rfloor} \pi(x^* - h_1 \tilde{u}, y^* - h_2 \tau_2 \tilde{v}) - \partial_y^{\lfloor \beta_2 \rfloor} \pi(x^* - h_1 \tilde{u}, y^*)\right)  \\
& \qquad + \frac{({-}h_1 \tilde{u})^{\lfloor \beta_1 \rfloor}}{\lfloor \beta_1 \rfloor!} \left(\partial_x^{\lfloor \beta_1 \rfloor} \pi(x^* - \tau_1 h_1 \tilde{u}, y^*) - \partial_x^{\lfloor \beta_1 \rfloor} \pi(x^*, y^*)\right) \Bigg) \dd \tilde{u} \dd \tilde{v}, \nonumber
\end{align}
where we once again use that \(K\) is a kernel of order \(M\).

Since \(\pi\) belongs to the H\"older space \(\mathcal{H}(\beta, \mathcal{L})\), we have the following bounds:

\[
|\partial_x^{\lfloor \beta_1 \rfloor} \pi(x^* - \tau_1 h_1 \tilde{u}, y^*) - \partial_x^{\lfloor \beta_1 \rfloor} \pi(x^*, y^*)| \leq \mathcal{L}_1 |\tau_1 h_1 \tilde{u}|^{\beta_1 - \lfloor \beta_1 \rfloor},
\]
\[
|\partial_y^{\lfloor \beta_2 \rfloor} \pi(x^* - h_1 \tilde{u}, y^* - h_2 \tau_2 \tilde{v}) - \partial_y^{\lfloor \beta_2 \rfloor} \pi(x^* - h_1 \tilde{u}, y^*)| \leq \mathcal{L}_2 |\tau_2 h_2 \tilde{v}|^{\beta_2 - \lfloor \beta_2 \rfloor}.
\]

Thus, \eqref{eq: order 11} is bounded above by

\[
\int_{\R^2} K(\tilde{u}) K(\tilde{v}) \left( \frac{\mathcal{L}_2(h_2 \tilde{v})^{\beta_2}}{\lfloor \beta_2 \rfloor!} + \frac{\mathcal{L}_1(h_1 \tilde{u})^{\beta_1}}{\lfloor \beta_1 \rfloor!}  \right) \dd \tilde{u} \dd \tilde{v} \leq C(h_2^{\beta_2} + h_1^{\beta_1}),
\]
which provides the desired result.
\end{proof}

We now turn to the proof of Theorem \ref{th: girsanov}, which serves as the central tool presented in Section \ref{s:prob}. While Section \ref{s:prob} outlines a heuristic argument for the result, here we provide a detailed and rigorous proof, filling in all the necessary steps.

\subsection{Proof of Theorem \ref{th: girsanov}}
\begin{proof}

Recall that, as discussed in the heuristic proof provided in Section \ref{s:prob}, the first step involves introducing the probability measure \(\mathbb{Q}\), under which the process \((N_u)_{0 \le u \le s}\) becomes a Poisson process with intensity \(\xi\). This transformation is achieved by introducing the random variable \(\mathbb{L}_s\) as defined in \eqref{eq: def L}. 

Next, we introduce \((\overline{\mathbb{L}}_s)\), a process closely related to \((\mathbb{L}_s)\), but where the dependence on \(\lambda_{T_1}\) is eliminated. Specifically, for \(u > T_1\), we know that 
\[
\lambda_u = \alpha \sum_{j \ge 2 : T_j {<} u}^q e^{-\beta (u - T_j)} + \xi + (\lambda_{T_1} + \alpha - \xi) e^{-\beta (u - T_1)}.
\]
Thus, for \(u > T_1\), we define
\begin{align}\label{eq : lambda bar sur}
\overline{\lambda}_u &= \alpha \sum_{j \ge 2 : T_j {<} u}^q e^{-\beta (u - T_j)} + \xi + (\lambda_0 + \alpha - \xi) e^{-\beta (u - T_1)}, 
\\ \label{eq : lambda bar sous}
\underline{\lambda}_u &= \alpha \sum_{j \ge 2 : T_j {<} u}^q e^{-\beta (u - T_j)} + \xi + \alpha e^{-\beta (u - T_1)}.
\end{align}
For \(u \le T_1\), we have \(\overline{\lambda}_u = \underline{\lambda}_u = \lambda_u\).
Since \(\lambda_{T_1} \in [\xi , \lambda_0 ]\), it follows that
\[
\underline{\lambda}_u \le \lambda_u \le \overline{\lambda}_u.
\]
Moreover, we immediately obtain the bound
\begin{equation}\label{eq: comparaison les lambda}
|\overline{\lambda}_u - \underline{\lambda}_u| \le (\lambda_0 - \xi) e^{-\beta (u - T_1)} \quad \text{for all } u > 0.    
\end{equation}
We define an upper bound for the process \((\mathbb{L}_s)_s\) by introducing \(\overline{\mathbb{L}}_s\) as in \eqref{eq: def overline L}, which we recall here for convenience:
\[
\overline{\mathbb{L}}_s = \left( \prod_{0 \le T_j \le s} \frac{\overline{\lambda}_{T_j}}{\xi} \right) \times \exp\left( - \int_0^s (\underline{\lambda}_u - \xi) \dd u \right).
\]
Clearly, \(\mathbb{L}_s \leq \overline{\mathbb{L}}_s\) almost surely for all \(s \geq 0\).

Let us now recall that, due to the change of measure we introduced, we have
\[
A^{(q)} = \mathbb{E}_{\mathbb{Q}} \left[ f^{(q)}(s - T_q, \ldots, T_2 - T_1, T_1) \, \mathbb{L}_s \, \one_{N_s = q} \mid \lambda_0 = y_0 \right].
\]
However, it will be more convenient going forward to use \(\mathbb{L}_{T_1 + s}\) as the change of measure rather than \(\mathbb{L}_s\), since we intend to work with the shifted process \(\widetilde{N}_\cdot = N_{T_1 + \cdot} - N_{T_1}\). 

Since the random variable inside the expectation is both \(\mathcal{F}_s\)-measurable and \(\mathcal{F}_{T_1 + s}\)-measurable, we can express it as
\begin{equation} \label{eq : lemma tech first hat g}
A^{(q)} = \mathbb{E}_{\mathbb{Q}} \left[ f^{(q)}(s - T_q, \ldots, T_2 - T_1, T_1) \, \mathbb{L}_{T_1 + s} \, \one_{N_s = q} \mid \lambda_0 = y_0 \right].
\end{equation}

We now use the inequality
\[
\mathbb{L}_{T_1 + s} \le \overline{\mathbb{L}}_{T_1 + s},
\]
where
\[
\overline{\mathbb{L}}_{T_1 + s} = \left( \prod_{j: T_j \leq T_1 + s} \frac{\overline{\lambda}_{T_j}}{\xi} \right) \times \exp\left( - \int_0^{T_1 + s} (\underline{\lambda}_u - \xi) \dd u \right),
\]
which can be rewritten as
\begin{align} \label{eq: lien L bar L tilde}
\overline{\mathbb{L}}_{T_1 + s} &= \frac{\overline{\lambda}_{T_1}}{\xi} \times \exp\left( - \int_0^{T_1} (\underline{\lambda}_u - \xi) \dd u \right) \times \widetilde{\mathbb{L}}_s,
\end{align}
where
\[
\widetilde{\mathbb{L}}_s := \left( \prod_{j \geq 2 : T_j \leq T_1 + s} \frac{\overline{\lambda}_{T_j}}{\xi} \right) \times \exp\left( - \int_{T_1}^{T_1 + s} (\underline{\lambda}_u - \xi) \dd u \right),
\]
or equivalently,
\[
\widetilde{\mathbb{L}}_s = \left( \prod_{j \geq 1 : \widetilde{T}_j \leq s} \frac{\overline{\lambda}_{T_1 + \widetilde{T}_j }}{\xi} \right) \times \exp\left( - \int_0^s (\underline{\lambda}_{T_1 + u} - \xi) \dd u \right),
\]
where we recall that \(\widetilde{T}_j = T_{j+1} - T_1\).

From the definitions in \eqref{eq : lambda bar sur}--\eqref{eq : lambda bar sous}, we see that both \(\overline{\lambda}_{T_1 + u}\) and \(\underline{\lambda}_{T_1 + u}\) are measurable with respect to \(\widetilde{\mathcal{F}}_u := \sigma(\widetilde{N}_v,~ v \le u)\). Consequently, \(\widetilde{\mathbb{L}}_s\) is also \(\widetilde{\mathcal{F}}_s\)-measurable, as required.

Note that
\[
\frac{\overline{\lambda}_{T_1}}{\xi} \exp\left(-\int_0^{T_1} (\underline{\lambda}_u - \xi) \dd u \right) \le \frac{\overline{\lambda}_{T_1}}{\xi} \le \frac{\lambda_0}{\xi}.
\]
Thus, we have the bound
\begin{equation}{\label{eq: 54.5}}
 \mathbb{L}_{T_1 + s} \le \overline{\mathbb{L}}_{T_1 + s} \le \widetilde{\mathbb{L}}_s \frac{\lambda_0}{\xi}.   
\end{equation}

Recalling \eqref{eq : lemma tech first hat g}, we can deduce:
\[
A^{(q)} \leq \mathbb{E}_{\mathbb{Q}} \left[ f^{(q)}(s - T_q, \ldots, T_2 - T_1, T_1) \, \one_{N_s = q} \, \widetilde{\mathbb{L}}_s \mid \lambda_0 = y_0 \right] \frac{\lambda_0}{\xi}.
\]
Now, as anticipated in \eqref{eq: meas tilde F}, we express \(\{N_s = q\}\) as \(\{T_1 \leq s - \widetilde{T}_{q-1}\} \cap \{\widetilde{N}_s = q - 1\}\), and using \(\widetilde{T}_j = T_{j+1} - T_1\), we obtain:
\[
A^{(q)} \leq \mathbb{E}_{\mathbb{Q}} \left[ \mathbb{E}_{\mathbb{Q}} \left[ f_1(s - \widetilde{T}_{q-1} - T_1) f_2^{(q)}(\widetilde{T}_{q-1}, \ldots, \widetilde{T_1}, T_1) \one_{T_1 \leq s - \widetilde{T}_{q-1}} \mid \widetilde{\mathcal{F}}_s \right] \one_{\{\widetilde{N}_s = q - 1\}} \widetilde{\mathbb{L}}_s \mid \lambda_0 = y_0 \right] \frac{y_0}{\xi},
\]
where we used the fact that both \(\wt{\mathbb{L}}_s\) and \(\wt{N}_s\) are \(\wt{\mathcal{F}}_s\)-measurable.

In the inner conditional expectation, the random variables \((\wt{T}_j)_{j=1, \dots, q-1}\) are \(\widetilde{\mathcal{F}}_s\)-measurable, since \(\widetilde{T}_{q-1} \leq s\). Consequently, this conditional expectation reduces to an integral with respect to the law of \(T_1\) conditioned on \(\widetilde{\mathcal{F}}_s\). Under the probability \(\mathbb{Q}\), the process \(N\) is a Poisson process, so \(T_1\) is independent of \(\widetilde{\mathcal{F}}_s = \sigma(N_{T_1 + u} - N_{T_1}, 0 \leq u \leq s)\), and its conditional law is exponential with parameter \(\xi\).

We can thus bound the inner conditional expectation by
\[
\int_0^{s - \wt{T}_{q-1}} f_1(s - \widetilde{T}_{q-1} - u) f_2^{(q)}(\widetilde{T}_{q-1}, \ldots, \wt{T_1}, u) \xi e^{-\xi u} \dd u.
\]
Since we assumed that \(f_2^{(q)}\) is bounded and non-zero only for \(u \in J^{(q)}\) (for some interval \(J^{(q)} \subset \mathbb{R}{_+}\)), this quantity is bounded by
\[
C \int_{[0, s - \wt{T}_{q-1}] \cap J^{(q)}} f_1(s - \widetilde{T}_{q-1} - u) \xi e^{-\xi u} \dd u \leq C \int_{s - \wt{T}_{q-1} - |J^{(q)}|}^{s - \widetilde{T}_{q-1}} f_1(s - \wt{T}_{q-1} - u) \xi \dd u = C \int_0^{|J^{(q)}|} \xi f_1(v) \dd v,
\]
where \(|J^{(q)}|\) is the length of the interval \(J\), and we used the change of variable \(v := s - \widetilde{T}_{q-1} - u\).

Using this bound, and $\sup_{q \ge 1} |J^{(q)}| \le \overline{J}$ we conclude, as in \eqref{eq: end prob}:
\[
\sum_{q=1}^\infty A^{(q)} \leq c \left( \int_0^{\overline{J}} f_1(v) \dd v \right) \sum_{q=1}^\infty \mathbb{E}_{\mathbb{Q}} \left[ \overline{\mathbb{L}}_{T_1 + s} \one_{\{\widetilde{N}_s = q - 1\}} \mid \lambda_0 = y_0 \right] = c \left( \int_0^{\overline{J}} f_1(v) \dd v \right) \mathbb{E}_{\mathbb{Q}} \left[ \overline{\mathbb{L}}_{T_1 + s} \mid \lambda_0 = y_0 \right].
\]
The proof is completed by applying the bound on the conditional expectation from Lemma \ref{L: upper bound Girsanov L1}.
\end{proof}

As outlined above, the proof of Theorem \ref{th: girsanov} relies on two key lemmas, Lemma \ref{L: upper bound Girsanov L1} and Lemma \ref{L: extension Leblanc}, both of which we will prove in detail in the following sections.

\subsection{Proof of Lemma \ref{L: upper bound Girsanov L1}}
\begin{proof}
As 	$\mathbb{E}_{\mathbb{Q}} \left[ {\mathbb{L}}_{T_1+s}\right] = 	\mathbb{E}_{\mathbb{P}} \left[1 \right]=1$, the matter of the proof is to study the deviation between ${\mathbb{L}}_{T_1+s}$ and its upper bound
	$\overline{{\mathbb{L}}}_{T_1+s}$. 
	
	Comparing \eqref{eq: def L}, \eqref{eq: def overline L} and using $\lambda_u=\overline{\lambda}_u=\underline{\lambda}_u$ for $u \le T_1$, we deduce
	\begin{equation*}
		\frac{\overline{\mathbb{L}}_{T_1+s}}{\mathbb{L}_{T_1+s}}
		= \left(\prod_{j : j \ge 2, T_j \le T_1+s} \frac{\overline{\lambda}_{T_j}}{\lambda_{T_j}} \right) \times
		\exp\left(\int_{T_1}^{T_1+s} (\lambda_u - {\underline{\lambda}_u}) \dd u\right).
	\end{equation*}
From $\underline{\lambda}_u \le \lambda_u \le {\overline{\lambda}_u}$ and \eqref{eq: comparaison les lambda} we have $\int_{T_1}^{T_1+s} (\lambda_u - {\underline{\lambda}_u})\dd u \le \int_0^\infty(\lambda_0-\xi)e^{-\beta u} \dd u\le (\lambda_0-\xi)/\beta $.
We deduce
\begin{align}\nonumber
			\mathbb{E}_{\mathbb{Q}} \left[\overline{\mathbb{L}}_{T_1+s} \mid \lambda_0=y_0\right] &=  
	\mathbb{E}_{\mathbb{Q}} \left[{\mathbb{L}}_{T_1+s} 
	\frac{\overline{\mathbb{L}}_{T_1+s}}{\mathbb{L}_{T_1+s}}
	\mid \lambda_0=y_0\right] 
	\\\nonumber & = 
	\mathbb{E}_{\mathbb{P}} \left[ 
	\frac{\overline{\mathbb{L}}_{T_1+s}}{\mathbb{L}_{T_1+s}}
	\mid \lambda_0=y_0\right]
	\\ \label{eq: majo EQ L by prod} &\le e^{\frac{y_0-\xi}{\beta}}	\mathbb{E}_{\mathbb{P}} \left[ 
\left(\prod_{j : j \ge 2, T_j \le T_1+s} \frac{
{\overline{\lambda}_{T_j}}}{
{\lambda_{T_j}}} \right)
	\mid \lambda_0=y_0\right].
	\end{align}
We have 
\begin{align}\label{eq: control ratio lambda}
\frac{\overline{\lambda}_{T_j}}{\lambda_{T_j}}
\le 1 + \frac{\overline{\lambda}_{T_j}-{\lambda}_{T_j}}{{\lambda}_{T_j}} &\le 1 + \frac{(\lambda_0-\xi)e^{-\beta(T_j-T_1)}}{\lambda_{T_j}} 
\\ \label{eq : control ratio avec xi}
&\le 1 + \frac{(\lambda_0-\xi)e^{-\beta(T_j-T_1)}}{\xi} 
\end{align} where we used \eqref{eq: comparaison les lambda} in the first line and $\lambda_{T_j} \ge \xi$ in the second line.
Now, let us fix $\varepsilon > 0$ arbitrarily small and $\overline{T}_\varepsilon >0$ such that $\frac{(\lambda_0-\xi)e^{-\beta \overline{T}_\varepsilon}}{\xi}  \le \varepsilon$. This allows us to write 
\begin{align*}
	\prod_{j : j \ge 2, T_j \le T_1+s} \frac{\overline{\lambda}_{T_j}}{\lambda_{T_j}} 
&	\le 
	\left( \prod_{j : j \ge 2, T_j \le T_1+\overline{T}_\varepsilon} \frac{\overline{\lambda}_{T_j}}{\lambda_{T_j}} \right)
	\times
	\left( \prod_{j : j \ge 2,   T_1+\overline{T}_\varepsilon < T_j \le T_1 + s} \frac{\overline{\lambda}_{T_j}}{\lambda_{T_j}} \right)
	\\
&   \le 
\left( \prod_{j : j \ge 2, T_j \le T_1+\overline{T}_\varepsilon}   1 + \frac{(\lambda_0-\xi)}{\lambda_{T_j}} \right)
\times
\left( \prod_{j : j \ge 2,   T_1+\overline{T}_\varepsilon < T_j \le T_1 + s} 1+\varepsilon \right)
\\
&\le 
\left( \prod_{j : j \ge 2, T_j \le T_1+\overline{T}_\varepsilon}   1 + \frac{(\lambda_0-\xi)}{\lambda_{T_j}} \right)
\times
\left( 1+\varepsilon \right)^{ N_{T_1+s} -N_{T_1+\overline{T}_\varepsilon}},
\end{align*}
where we used \eqref{eq: control ratio lambda}--\eqref{eq : control ratio avec xi} in the second line.
Recalling the notation $\widetilde{N}_u=N_{T_1+u}-N_{T_1}$, and $\widetilde{N}((a,b])=\widetilde{N}_b-\widetilde{N}_a$ it is
\begin{align*}
	\prod_{j : j \ge 2, T_j \le T_1+s} \frac{\overline{\lambda}_{T_j}}{\lambda_{T_j}} 
	\le 
\left( \prod_{j : j \ge 2, T_j \le T_1+\overline{T}_\varepsilon}   1 + \frac{(\lambda_0-\xi)}{\lambda_{T_j}} \right)
	\times
	\left( 1+\varepsilon \right)^{ \widetilde{N}((\overline{T}_\varepsilon,s])}.
\end{align*}
Now, we use $\lambda_{T_j} \ge \lambda_{T_1} e^{-\beta(T_j-T_1)} + (1-e^{-\beta(T_j-T_1)})\xi +
\sum_{1 \le l <j} \alpha e^{-\beta (T_j-T_l)} \ge \xi + \alpha (j-1) e^{-\beta \overline{T}_\varepsilon}$ if $T_j \le T_1+\overline{T}_\varepsilon$, and deduce
\begin{align*}
	\prod_{j : j \ge 2, T_j \le T_1+s} \frac{\overline{\lambda}_{T_j}}{\lambda_{T_j}} 
&	\le 
	\left( \prod_{j : j \ge 2, T_j \le T_1+\overline{T}_\varepsilon}   1 + \frac{(\lambda_0-\xi)}{\xi + \alpha (j-1) e^{-\beta \overline{T}_\varepsilon}} \right)
	\times
	\left( 1+\varepsilon \right)^{ \widetilde{N}((\overline{T}_\varepsilon,s])}
	\\
&	\le C \left( \prod_{j =1}^{\widetilde{N}[0,\overline{T}_\varepsilon]} 1 + \frac{C}{j} \right)
	\times
	\left( 1+\varepsilon \right)^{ \widetilde{N}((\overline{T}_\varepsilon,s])}
\end{align*}
for some constant $C>0$. Then, we use the simple inequality $\prod_{j =1}^{q}(1+\frac{C}{{j}}) \le C' q^{C'}$ for some constant $C'$ independent of $q \ge 1$. It entails,
\begin{equation*}
	\prod_{j : j \ge 2, T_j \le T_1+s} \frac{\overline{\lambda}_{T_j}}{\lambda_{T_j}} 
	\le C' (\widetilde{N}[0,\overline{T}_\varepsilon])^{C'}
	\times
	\left( 1+\varepsilon \right)^{ \widetilde{N}((\overline{T}_\varepsilon,s])}.
\end{equation*}
	Hence, using \eqref{eq: majo EQ L by prod},
	\begin{align*}
		\mathbb{E}_{\mathbb{Q}} \left[\overline{\mathbb{L}}_{T_1+s} \mid \lambda_0=y_0\right]
		&\le C' e^{\frac{y_0-\xi}{\beta}}
				\mathbb{E}_{\mathbb{P}} \left[
			(\widetilde{N}[0,\overline{T}_\varepsilon])^{C'}
			\times
			\left( 1+\varepsilon \right)^{ \widetilde{N}((\overline{T}_\varepsilon,s])}		
				\mid \lambda_0=y_0\right]
					\\	&\le C'e^{\frac{y_0-\xi}{\beta}}
				\mathbb{E}_{\mathbb{P}} \left[
			(\widetilde{N}[0,\overline{T}_\varepsilon])^{2C'} \mid \lambda_0=y_0\right]^{1/2}
				\times
				\mathbb{E}_{\mathbb{P}} \left[	\left( 1+\varepsilon \right)^{ 2\widetilde{N}((\overline{T}_\varepsilon,s])}				
				\mid \lambda_0=y_0\right]^{1/2},
	\end{align*}
	where we used the Cauchy-Schwarz inequality in the last line.
Using the Markovian property of the exponential Hawkes process, we know that  $\widetilde{N}$ is a Hawkes process with starting intensity $\widetilde{\lambda}_0=\lambda_{T_1}+\alpha=\xi+(\lambda_0-\xi)e^{-\beta T_1}+\alpha\le \lambda_0 + \alpha$. 
If $\varepsilon$ is small enough, we can use Lemma \ref{L: extension Leblanc} on the exponential moments of Hawkes process with $K=(1+\varepsilon)^2$. As already mentioned, this lemma is an extension to the non stationary case of the result in \cite{leblanc2024exponential}. It yields :
\begin{multline*}
		\mathbb{E}_{\mathbb{Q}} \left[\overline{\mathbb{L}}_{T_1+s} \mid \lambda_0=y_0\right]
		\le C' e^{\frac{y_0-\xi}{\beta}}
		\mathbb{E}_{\mathbb{P}} \left[
		(\widetilde{N}[0,\overline{T}_\varepsilon])^{2C'} \mid \lambda_0=y_0\right]^{1/2}
		\\\times
	\exp\left( \frac{1}{2} \big( (1+\varepsilon)^{2(1+\frac{2}{1-\beta/\alpha})}-1 \big) (s-\overline{T}_\varepsilon)(y_0+\alpha) \right). 
\end{multline*}
From Lemma \ref{L: extension Leblanc}, we see that the Hawkes process admits finite moments of any order, and then by setting $q(\varepsilon)= \frac{1}{2}(1+\varepsilon)^{2(1+\frac{2}{1-\beta/\alpha})}-1 $, we have
\begin{equation*}
	\mathbb{E}_{\mathbb{Q}} \left[\overline{\mathbb{L}}_{T_1+s} \mid \lambda_0=y_0\right]
	\le {C_\varepsilon} \exp\left( q(\varepsilon) (s-\overline{T}_\varepsilon)(y_0+\alpha) \right) \le {C_\varepsilon} 	\exp\left( q(\varepsilon)s (y_0+\alpha) \right),
\end{equation*}
for some constant {$C_\varepsilon= C' e^{\frac{y_0-\xi}{\beta}}
		\mathbb{E}_{\mathbb{P}} \left[
		(\widetilde{N}[0,\overline{T}_\varepsilon])^{2C'} \mid \lambda_0=y_0\right]^{1/2} <\infty$. Observe that this constant is bounded, on compact sets, by a constant that does not depend on $y_0$. This is quite straightforward to be remarked, thanks to the statement of Lemma \ref{L: extension Leblanc}.}	

{Let $\varepsilon'>0$. 
As $\varepsilon$ is arbitrarily small and $q(\varepsilon)\xrightarrow{\varepsilon\to 0}0$, 
we can fix $\varepsilon$ such that $q(\varepsilon)\le \varepsilon'$.
We have proved}
\begin{equation*}
	\mathbb{E}_{\mathbb{Q}} \left[\overline{\mathbb{L}}_{T_1+s} \mid \lambda_0=y_0\right]
	\le C_\varepsilon 	\exp\left( {\varepsilon'} s \right),
\end{equation*}
	 for any fixed {$\varepsilon'>0$,} which is the lemma. 
\end{proof}

\subsection{Proof of Lemma \ref{L: extension Leblanc}}
\begin{proof}
In \cite{leblanc2024exponential}, it is shown that if $\overline{N}$ is a stationary exponential Hawkes process with parameters $(\overline{\xi}, \alpha, \beta)$, where $\alpha/\beta < 1$, the following inequality holds:

\begin{equation}{\label{eq: res Leblanc}}
 \E\left[K^{\overline{N}[0,t]} \right] \le \exp\left( \left( K^{1+\frac{2}{1-\alpha/\beta}} - 1 \right) t \overline{\xi} \right),   
\end{equation}

for $1 < K \le \left(\frac{\beta}{2\alpha} + \frac{1}{2}\right)^{1+\frac{2}{1-\alpha/\beta}}$. 

To extend this result to a non-stationary version of the Hawkes process, we construct a coupling between two Hawkes processes using the cluster representation of the process, as described, for instance, in \cite{hawkesClusterProcessRepresentation1974} or \cite{mollerPerfectSimulationHawkes2005c}.

We can construct the non-stationary Hawkes process $\wt{N}$ on $[0, t]$ by first generating arrival times of an inhomogeneous Poisson process $\wt{P}$ with intensity:

\[\tilde{\mu}_s = \left( y_0 e^{-\beta s} + \xi (1 - e^{-\beta s}) \right) \one_{s \geq 0}.
\]

These arrival times, denoted by $T_1^{(0)}, \dots, T_m^{(0)}$, are referred to as "immigrants". Each immigrant $T_j^{(0)}$ generates a cluster of points through successive generations of offspring, following Poisson processes with intensity $s \mapsto {\alpha} e^{-\beta(s-d)} \one_{s > d}$, where $d$ is the birth date of the ancestor. These clusters are independent for each immigrant. The Hawkes process is the union of the immigrants and all cluster points.

A stationary version $\overline{N}$ of the process, with parameters $(\overline{\xi}, \alpha, \beta)$, can be obtained in the same way, except the immigration layer is drawn from a homogeneous Poisson process $\overline{P}$ with constant intensity $\overline{\xi}$ over the entire real line $\R$. If $y_0 \geq \xi$, we choose $\overline{\xi} = y_0$; otherwise, we set $\overline{\xi} = \xi$. In any case, the intensity of the homogeneous Poisson process $\overline{P}$ dominates that of the inhomogeneous process $\wt{P}$. 

By applying a thinning procedure (see Theorem 1 in \cite{lewisSimulationNonhomogeneousPoisson1979}), it is possible to construct a coupling of $\overline{P}$ and $\widetilde{P}$ such that all jump times of $\widetilde{P}$ are included in the jump times of $\overline{P}$. By using the same random clusters for the immigrants common to both immigration processes, we find that the jump times of the non-homogeneous Hawkes process are embedded within those of the stationary Hawkes process. Thus, we have $\wt{N}[0,t] \leq \overline{N}[0,t]$. The lemma then follows immediately from \eqref{eq: res Leblanc} with $\overline{\xi} = y_0$.    
\end{proof}

We now proceed to the proof of the final two lemmas, whose statements are presented in Section \ref{sec:proofs}, as they have proven essential for establishing the variance bound in the case where $y^* > \xi$. Hence, to conclude the paper, we provide the proofs of Lemmas \ref{L: g pour q zero} and \ref{L: g smaller g hat}.

\subsection{Proof of Lemma \ref{L: g pour q zero}}
\begin{proof}
Under the assumption that $N_s=0$, the process $\lambda$ is deterministic on $[0,s]$ given by 
	$\lambda_t= \xi+(\lambda_0-\xi)e^{-\beta t}$. Since the support of $\mathbb{K}(y^*-\cdot)$ is the interval $[y^*-h_2,y^*+h_2]$, we deduce that
	$g^{(0)}(s)=0$ as soon as $\lambda_s < y^*-h_2$.
	This is equivalent to 
	\begin{equation*}
		e^{-\beta s}< \frac{y^* -h_2-\xi}{\lambda_0-\xi}, 
	\end{equation*}
	which gives using $\lambda_0=y_0$,
	\begin{equation*}
		s > -\frac{1}{\beta} \ln \left( \frac{y^* -h_2-\xi}{y_0-\xi}\right)
	= -\frac{1}{\beta} \ln \left(1 + \frac{y^* - y_0 - h_2}{y_0 - \xi}\right). 
	\end{equation*}
	Now, we use that $y^* - y_0 - h_2 \in [-2h_2,0]$ and $y_0 - \xi = y_0-y^*+y^* - \xi \ge -h_2+5h_2=4h_2$ to deduce that $ \frac{y^* - y_0 - h_2}{y_0 - \xi} \in [-1/2,0]$. As $-\ln (1+u) \le - 2 u$ for $u \in [-1/2,0]$, we can write that $g^{(0)}(s)=0$ as soon as 
	\begin{equation*}
		s > -\frac{2}{\beta} \times \frac{y^* - y_0 - h_2}{y_0-\xi}
	\end{equation*}
	and thus as soon as $s \ge \frac{4}{\beta} \frac{h_2}{y_0-\xi}$ using again $|y^* -h_2-\xi| \le 2h_2$. As under our assumption we have $y_0-\xi=(y^*-\xi)(1+\frac{y_0-y^*}{y^*-\xi})
	\ge (y^*-\xi)\frac{4}{5}$, we see that a sufficient condition for $g^{(0)}(s)=0$ is  $s \ge  \frac{5}{\beta}\frac{h_2}{y^*-\xi}$. Observe that, in the above, we are not seeking a lower bound on \( s \). Instead, we are identifying a region of \( \mathbb{R} \) where a certain condition holds, specifically \( g^{(0)}(s) = 0 \). In fact, if we know that \( g^{(0)}(s) = 0 \) for any \( s > a \) and we find some \( b > a \), then it is clear that \( g^{(0)}(s) = 0 \) also holds for \( s > b \) (since \( b \) is greater than \( a \)).
	The lemma is proved.
    
\end{proof}

\subsection{Proof of Lemma \ref{L: g smaller g hat}}
\begin{proof}
We know that on $\{N_s=q\}$, $\lambda_s=\alpha \sum_{i=1}^q  e^{-\beta (s-T_i)} +\xi+(\lambda_0-\xi)e^{-\beta s}$ {almost surely}. We use that $\mathbb{K}_{h_2}(y^*-\lambda_s) \neq 0$ only if $\lambda_s \in [y^*-h_2,y^*+h_2]$.
		Thus $\mathbb{K}_{h_2}(y^*-\lambda_s) \neq 0$ implies
		$\alpha \sum_{i=1}^q e^{-\beta (s-T_i)} +\xi+(\lambda_0-\xi)e^{-\beta s} \in [y^*-h_2,y^*+h_2]$, which is 
		$$ e^{-\beta s} \left[ \alpha e^{\beta T_1} \left(1+\sum_{j=2}^q e^{\beta (T_j-T_1)}\right) + \lambda_0-\xi \right] \in [\ell_g,\ell_d].
		$$
		It yields $ |\mathbb{K}_{h_2}(y^*-\lambda_s)|	\le \frac{\norm{{K}}_\infty}{h_2} \one_{\{e^{\beta T_1} \left(1+\sum_{j=2}^q e^{\beta (T_j-T_1)}\right) \in I \}}$ with $I=[e^{\beta s} \ell_g + \xi -\lambda_0, e^{\beta s} \ell_d + \xi -\lambda_0]$. The result follows, recalling that $L_s=T_q$ is the last jump occurring before time $s$ on $\{N_s=q\}$.
\end{proof}

\end{document}